%% file: main.tex
\documentclass{article}

\usepackage{amsmath,amssymb,amsthm}
\usepackage{epsfig}
\usepackage[utf8x]{inputenc}
\usepackage{psfrag,textcomp}
\usepackage{algorithm,algo_olm}
\usepackage{color}
\usepackage{mathtools}
\usepackage{graphicx,import}
\usepackage{pgfplots}
\usepackage{caption} 
\usepackage{subcaption}
\usepackage{tikz}
\usepackage{a4wide}

\newtheorem{theorem}{Theorem}
\newtheorem{proposition}[theorem]{Proposition}
\newtheorem{remark}[theorem]{Remark}
\newtheorem{lemma}[theorem]{Lemma}
\newtheorem{coro}[theorem]{Corollary}
\newtheorem{property}[theorem]{Property}

\newtheorem{interm_result}[theorem]{Intermediate result}

\pgfplotsset{compat=newest}

\newcommand{\powerset}{\raisebox{.15\baselineskip}{\Large\ensuremath{\wp}}}

\newcommand{\Kappa}{\overline{\kappa}}

\newcommand\blfootnote[1]{%
  \begingroup
  \renewcommand\thefootnote{}\footnote{#1}%
  \addtocounter{footnote}{-1}%
  \endgroup
}

\begin{document}

\title{Nonintrusive approximation of parametrized limits of matrix power algorithms -- application to matrix inverses and log-determinants}
\author{Fabien Casenave${}^{*}$, Nissrine Akkari${}^{*}$, Alexandre Charles${}^{*}$, Christian Rey${}^{*}$\blfootnote{Safran Tech, Rue des Jeunes Bois, Ch\^{a}teaufort CS 80112 - 78772 Magny-les-Hameaux Cedex - France}}
\date{\today}
\maketitle
\begin{abstract}  
We consider in this work quantities that can be obtained as limits of powers of parametrized matrices, for instance the inverse matrix or the logarithm of the determinant. Under the assumption of affine dependence in the parameters, we use the Empirical Interpolation Method (EIM) to derive an approximation for powers of these matrices, from which we derive a nonintrusive approximation for the aforementioned limits. We derive upper bounds of the error made by the obtained formula. Finally, numerical comparisons with classical intrusive and nonintrusive approximation techniques are provided: in the considered test-cases, our algorithm performs well compared to the nonintrusive ones.
\end{abstract}

\maketitle

\section{Introduction}

Many models in physics, biology or engineering involve partial differential equations, which are nowadays mainly solved numerically. In many cases, a single solution is not enough, as we are interested in the behavior of the solution when some chosen parameters vary. For instance, in sensitivity analyses, optimization or uncertainty quantification, the solution has to be computed a large number of times. In this many-queries context, model order reduction techniques have proved to allow large improvements in computational costs.

A number of techniques and methods can be grouped under the heading of Reduced Order Models (ROM). First, one can simply consider taking a coarser mesh, or making use of symmetries in the problem. Then, one can use methods from the machine learning community, where a meta-model is constructed as an interpolation or regression of the solutions or quantities of interest over the parameter set.
These techniques are nonintrusive since they use the numerical solver as a black-box, see~\cite{reviewmachinelearning2,reviewmachinelearning3,reviewmachinelearning4,reviewmachinelearning1,reviewmachinelearning5} for reviews of machine learning regression methods. Finally, a third class of ROM consists in solving the partial differential equation (or an approximation of it) on a small dimensional subspace, so that the computational cost of solving the reduced model is orders of magnitude smaller than that of the full-scale model. For instance, one can use the Proper Orthogonal Decomposition (POD)~\cite{sirovich,POD} or the Reduced-Basis method~\cite{prud1,patera,Machiels,RB,prud2,sen,naturalnorm,3,2,yano}. The Proper Generalized Decomposition (PGD) usually expresses the solution as a function of space and time, and the parameters of the model, seen here as coordinates. This function is approximated as a sum of tensor products, see~\cite{PGD, PGD2, PGD3, PGD4}. These methods are generally very intrusive to the considered computational code, since they need to modify the assembly routines of the operators. Efforts have been spent to mitigate these intrusivity requirements~\cite{casenave2015nonintrusive, chakir:hal-00855906, FLD:FLD4263}, but we still need to manipulate at least the matrices or meshes.

In this work, we consider a family of parametrized invertible matrices, such that the parameter dependence is affine, and we are interested in the nonintrusive approximation of quantities obtained as limits of powers of parametrized matrices, for instance the inverse matrices or the logarithm of the determinant (log-det), which we express as linear combinations of these inverses or log-det computed at given parameter values. 
For instance, many evaluations of the log-det of a positive-definite matrix are required for maximum likelihood estimation in Gaussian process regression, see~\cite{logdet1,logdet2}.
The proposed algorithm computes the coefficients of these linear combinations efficiently (namely in a computational complexity independent of the size of the matrices), and is nonintrusive, in the sense that it resorts only to the evaluation of the quantities of interests, namely the inverse or the log-det, like most machine learning methods.
This is an offline/online procedure.  A computationally demanding stage is first carried out, the offline stage, where the high-fidelity model is solved a certain number of times and some information on the parameter dependence of the model is learned. This information is then exploited in the online stage, in a computationally cheap fashion, where the approximation is computed rapidly and potentially for a large number of parameter values.

In Section~\ref{sec:nonintrusive_EIM} is proposed an interpolation formula of the power of  parametrized matrices based on the Empirical Interpolation Method (EIM)~\cite{Barrault,Maday}.
Then in Section~\ref{sec:PowerAlgo}, nonintrusive approximations of limits to certain power algorithms are obtained using the interpolation of power of matrices, namely the inverse and the log-det of parametrized matrices. In Section~\ref{sec:convergence}, upper bounds of the error made by these approximations are derived. Finally, in Section~\ref{sec:numerics} are presented numerical comparisons between the proposed nonintrusive approximation algorithm and classical intrusive and nonintrusive methods.

\section{Approximation of powers of parametrized matrices}
\label{sec:nonintrusive_EIM}

Let $d\in\mathbb{N}$ and $\mu\in\mathcal{P}$ be a parameter, where the parameter set $\mathcal{P}$ is a compact subset of $\mathbb{R}^r$, $r\in\mathbb{N}^*$. Consider $\{A_\mu\}_{\mu\in\mathcal{P}} \subset \mathbb{R}^{\mathcal{N}\times \mathcal{N}}$ a set of parametrized square matrices and assume the following affine decomposition for each element of the set:
\begin{equation}
\label{eq:decompA0}
A_\mu=\sum_{l=1}^d\alpha_l(\mu)A_l, ~\mu\in\mathcal{P}, 
\end{equation}
\noindent where we suppose that the family of square matrices $\left\{A_l\right\}_{1\leq l\leq d}$ is independent of $\mu$.
Hence, the matrix $A_\mu$ depends on $\mu$ only through the coefficients $\alpha_l:\mathcal{P}\rightarrow \mathbb{R}$.

Let $m\in\mathbb{N}$. We propose to derive an offline/online procedure to compute an approximation of $A_\mu^p$ for $1\leq p\leq m$ and $\mu\in\mathcal{P}$ of the following form:
\begin{equation}
\label{eq:online_sol}
A_\mu^p\approx \sum_{l=1}^t \lambda_l(\mu) A_{\mu_l}^{p}, ~\mu\in\mathcal{P},
\end{equation}
where $\left\{\mu_l \right\}_{1\leq l\leq t}$ is determined during the offline stage whereas the applications $\lambda_l:\mathcal{P}\rightarrow\mathbb{R}$ are computed during the online stage. As we shall see later, this expression will be used to obtain efficient approximations for the inverse and the log-det of $A_\mu$.

Consider the decomposition~\eqref{eq:decompA0} and take the $p$-th power of the equation. In the general case, the matrices $A_i$, $1\leq i\leq d$ do not commute, which prevents us from the use of the multinomial theorem. Thus,
\begin{equation}
\label{eq:grossesomme}
A_\mu^p=\sum_{l_1=1}^d\sum_{l_2=1}^d\cdots\sum_{l_p=1}^d\left(\prod_{i=1}^p\alpha_{l_i}(\mu)\right)\left(\prod_{i=1}^pA_{l_i}\right),~\mu\in\mathcal{P}.
\end{equation}
In the following, we factorize the sum according to the products of $A_l$ matrices. Let $t\in \mathbb{N}$. Denote a multi-index $\vec{s}=\left(s_1,s_2,\cdots,s_t\right)\in\mathbb{N}^t$ and define its weight $|\vec{s}|:=\sum_{l=1}^t s_l$. Finally, denote $\Kappa_{m,d}=\left\{\vec{k}\in\textnormal{\textlbrackdbl}0;m\textnormal{\textrbrackdbl}^d\textnormal{ such that }|\vec{k}|\leq m\right\}$.

\begin{lemma}
\label{lemmaconstructif}
Let $p,\,d,\,m\in\mathbb{N}$ and $\{A_\mu\}_{\mu\in\mathcal{P}} \subset \mathbb{R}^{\mathcal{N}\times \mathcal{N}}$ a set of parametrized matrices satisfying~\eqref{eq:decompA0}. There exists $\{T_{\vec{k},p}\}_{\vec{k}\in\Kappa_{m,d},0\leq p\leq m}\in\mathbb{R}^{\mathcal{N}\times \mathcal{N}}$, independent of $\mu$, such that the following equality holds:
\begin{equation}
\label{eq:decompA2s}
A_\mu^p=\sum_{\vec{k}\in\Kappa_{m,d}}g(\vec{k},\mu) {T}_{\vec{k},p},~\mu\in\mathcal{P},~1\leq p\leq m,
\end{equation}
where
\begin{equation}
\label{eq:formula_g}
g(\vec{k},\mu)=\prod_{l=1}^d\alpha_l^{k_l}(\mu).
\end{equation}
\end{lemma}
\renewcommand*{\proofname}{\textbf{Proof of Lemma~\ref{lemmaconstructif}}}
We clarify here that Lemma~\ref{lemmaconstructif} contains an existence result for the matrices ${T}_{\vec{k},p}$, that do not need to be computed for the method to be carried out in practice. Equation~\eqref{eq:decompA2s} indicates that the function $g\mapsto A^p$ is linear, and the idea is to use an EIM approximation of $g(\vec{k},\mu)$ to readily obtain an approximation of $A^p_\mu$.
Working on $g(\vec{k},\mu)$ instead of the matrix coefficients $\left(A^p_\mu\right)_{i,j}$, $1\leq i,j\leq\mathcal{N}$, will enable us to construct nontrusive approximations for the inverse or log-det of $A^p_\mu$, as we will see in Section~\ref{sec:PowerAlgo}.

\begin{proof}
Define the family of applications:
\begin{equation}
\mathcal{F}_{p,d}:\quad\left\lbrace
\begin{array}{lll}
\textnormal{\textlbrackdbl}1;d\textnormal{\textrbrackdbl}^p&\rightarrow&\textnormal{\textlbrackdbl}0;p\textnormal{\textrbrackdbl}^d\\
\vec{s}&\mapsto&\left(\#\left\{s_i=j,~1\leq i\leq p\right\}\right)_{j=1,\cdots,d}
\end{array}\right.\qquad\forall d,q\in \mathbb{N}
\end{equation}
By construction, $|\mathcal{F}_{p,d}(\vec{s})|=p$ for all $\vec{s}\in \textnormal{\textlbrackdbl}1;d\textnormal{\textrbrackdbl}^p$.
Take $d=2$, $p=3$ for example. Then $\mathcal{F}_{3,2}((2,2,2)) = (0,3)$ and $\mathcal{F}_{3,2}((1,2,2))=\mathcal{F}_{3,2}((2,1,2))=\mathcal{F}_{3,2}((2,2,1))=(1,2)$.

Define now the reciprocal applications:
\begin{equation}
\mathcal{I}_{p,d}:\quad
\left\lbrace
\begin{array}{lll}
\textnormal{\textlbrackdbl}0;p\textnormal{\textrbrackdbl}^d&\rightarrow&\powerset\left(\textnormal{\textlbrackdbl} 1 ; d \textnormal{\textrbrackdbl}^p\right)\\
\vec{k}&\mapsto&\left\lbrace  \vec{s}\in\textnormal{\textlbrackdbl}1;d\textnormal{\textrbrackdbl}^p\quad :\quad \mathcal{F}_{p,d}\left(\vec{s}\right)=\vec{k} \right\rbrace
\end{array}
\right.\qquad\forall d,q\in \mathbb{N}
\end{equation}
where $\powerset\left(\textnormal{\textlbrackdbl} 1 ; d \textnormal{\textrbrackdbl}^p\right)$ denotes the power set of $\textnormal{\textlbrackdbl} 1 ; d \textnormal{\textrbrackdbl}^p$.
Take $d=2$, $p=3$ for example. Then, $\mathcal{I}_{3,2}((0,3))=\left\{(2,2,2)\right\}$ or $\mathcal{I}_{3,2}((1,2))=\left\{(1,2,2),(2,1,2),(2,2,1)\right\}$.

Using the introduced notation, Equation~\eqref{eq:grossesomme} can be reordered in the following form:
\begin{equation}
A_\mu^p=\sum_{\vec{k}\in\textnormal{\textlbrackdbl}0;p\textnormal{\textrbrackdbl}^d:|\vec{k}|=p}\left(\prod_{l=1}^d\alpha_l^{k_l}(\mu)\right)\sum_{\vec{s}\in \mathcal{I}_{p,d}(\vec{k})}\prod_{i=1}^p{A_{s_i}},\quad\forall\mu\in\mathcal{P}.
\end{equation}

Notice that if the matrices $A_l$, $1\leq l\leq d$, were commuting, we could have simply applied the multinomial theorem to get
\begin{equation}
A_\mu^p=\sum_{\vec{k}\in\textnormal{\textlbrackdbl}0;p\textnormal{\textrbrackdbl}^d:|\vec{k}|=p}\frac{p!}{\prod_{l=1}^dk_l!}\left(\prod_{l=1}^d\alpha_l^{k_l}(\mu)\right)\left(\prod_{l'=1}^d A_{l'}^{k_{l'}}(\mu)\right),\quad\forall\mu\in\mathcal{P}.
\end{equation}

Recall the notation $g(\vec{k},\mu)=\prod_{l=1}^d\alpha_l^{k_l}(\mu)$ and 
denote 
\begin{equation}
\label{eq:formuleTtilde}
\tilde{T}_{\vec{k},p}:=\sum_{\vec{s}\in\mathcal{I}_{p,d}(\vec{k})}\prod_{i=1}^p{A_{s_i}},
\end{equation}
so that 
\begin{equation}
\label{eq:decompA}
A_\mu^p=\sum_{\vec{k}\in\textnormal{\textlbrackdbl}0;p\textnormal{\textrbrackdbl}^d:|\vec{k}|=p}g(\vec{k},\mu) \tilde{T}_{\vec{k},p},~\mu\in\mathcal{P}.
\end{equation}
Denote now for a general $\vec{k}\in \textnormal{\textlbrackdbl}0;p\textnormal{\textrbrackdbl}^d$ (not restricted to only $|\vec{k}|=p$):
\begin{equation}
\label{eq:formulaTkp}
{T}_{\vec{k},p}:=\left\{
\begin{aligned}
\tilde{T}_{\vec{k},p}\textnormal{ if }|\vec{k}|=p,\\
0\textnormal{ otherwise.}
\end{aligned}
\right.
\end{equation}
Let $m\in\mathbb{N}$. The $p$-exponent in~\eqref{eq:decompA} can be parametrized using
\begin{equation}
\label{eq:decompA2ss}
A_\mu^p=\sum_{\vec{k}\in\Kappa_{m,d}}g(\vec{k},\mu) {T}_{\vec{k},p},~\mu\in\mathcal{P},\quad 1\leq p\leq m,
\end{equation}
where we recall that $\Kappa_{m,d}=\left\{\vec{k}\in\textnormal{\textlbrackdbl}0;m\textnormal{\textrbrackdbl}^d\textnormal{ such that }|\vec{k}|\leq m\right\}$, which concludes the proof.
\end{proof}

To illustrate Lemma~\ref{lemmaconstructif}, consider the case $p=2$ and $d=2$. In this case, $\left\{\vec{k}\in\textnormal{\textlbrackdbl}0;2\textnormal{\textrbrackdbl}^2:|\vec{k}|=2\right\}=\left\{(1,1),(0,2),(2,0)\right\}$, and from Equations~\eqref{eq:formulaTkp}-\eqref{eq:decompA2ss}, there holds
\begin{equation}
\label{eq:exampleLemma1}
A_\mu^2=g((1,1),\mu) {T}_{(1,1),2}+g((0,2),\mu) {T}_{(0,2),2}+g((2,0),\mu) {T}_{(2,0),2}. 
\end{equation}
Since $\mathcal{F}_{2,2}\left((1,2)\right)=\mathcal{F}_{2,2}\left((2,1)\right)=(1,1)$, there holds $\mathcal{I}_{2,2}\left((1,1)\right)=\left\{(1,2),(2,1)\right\}$, and from Equation~\eqref{eq:formuleTtilde}, $\tilde{T}_{(1,1),2}=A_1A_2+A_2A_1$. In the same fashion, we compute $\mathcal{F}_{2,2}\left((2,2)\right)=(0,2)$ leading to $\mathcal{I}_{2,2}\left((0,2)\right)=\left\{(2,2)\right\}$ and $\tilde{T}_{(0,2),2}=A_2^2$, as well as $\mathcal{F}_{2,2}\left((1,1)\right)=(2,0)$ leading to $\mathcal{I}_{2,2}\left((2,0)\right)=\left\{(1,1)\right\}$ and $\tilde{T}_{(2,0),2}=A_1^2$. Using the formula~\eqref{eq:formula_g} for $g$ in~\eqref{eq:exampleLemma1} leads to the 
known expression $A_{\mu}^2=\alpha_1(\mu)\alpha_2(\mu)\left(A_1A_2+A_2A_1\right)+\alpha_1(\mu)^2A_1^2+\alpha_2(\mu)^2A_2^2$. 

It is known that $\#\{\vec{k}\in\textnormal{\textlbrackdbl}0;p\textnormal{\textrbrackdbl}^d\textnormal{ such that }|\vec{k}|=p\}=\frac{(p+d-1)!}{(d-1)!p!}$. Then, the number of terms in~\eqref{eq:decompA2s}, namely 
$Q_{m,d}:=\#\Kappa_{m,d}$, equals $\frac{1}{(d-1)!}\sum_{k=0}^m\frac{(k+d-1)!}{k!}$. Notice that for $d\geq 2$, $Q_{m,d}\leq \frac{m}{\left(d-1\right)!}\left(m+1\right)\left(m+2\right)\cdots\left(m+d-1\right)=P_d(m)$, where $P_d(m)$ is a polynomial of degree $d$ in $m$.

As discussed earlier, we
carry out the Empirical Interpolation Method (EIM) on the $g(\vec{k},\mu)$, see Algorithm~\ref{algoEIM} for a description of the offline stage of EIM on this function. In Algorithm~\ref{algoEIM}, $\delta^l g:=I^l(g)-g$, with $I^l(g)$ denoting the rank-$l$ EIM approximation, defined by
\begin{equation}
\label{eq:EIM0}
I^l(g)(\vec{k},\mu)=\sum_{l'=1}^{l}\beta^l_{l'}(\mu)q^{l'}(\vec{k}),
\end{equation}
where $\beta^l(\mu)$ solves
\begin{equation}
\label{eq:linear_sys}
\sum_{l''=1}^{l} B_{l',l''}^{l}\beta^{l}_{l''}(\mu)=g(\vec{k}_{l'},\mu),\qquad 1\leq l'\leq l.
\end{equation}
The quantities $B^l\in\mathbb{R}^{l\times l}$, $q^l:\mathcal{P}\to\mathbb{R}^{l}$, $\vec{k}_l\in\Kappa_{m,d}$, $\mu_l\in\mathcal{P}_{\rm sample}$, for all $1\leq l\leq N^{\rm EIM}$, are constructed during the offline stage in Algorithm~\ref{algoEIM}, where $N^{\rm EIM}$ is the number of terms selected by the EIM. In practice, $N^{\rm EIM}$ is not \textit{a priori} specified, but results from a stopping criterion on the maximum current error $(\delta^l g)(\vec{k}_{l+1},\mu_{l+1})$ made by the approximation.

\begin{algorithm}[h!]
	\caption{Offline stage of the EIM}
	\label{algoEIM}
	\begin{algorithmic}[1]
	\STATE {Choose a fine finite set $\mathcal{P}_{\rm sample}\subset\mathcal{P}$}
	\STATE {Set $l:=1$}
        \STATE {Compute $\displaystyle \mu_1:=\underset{\mu\in \mathcal{P}_{\text{sample}}}{\textnormal{argmax}}\|g(\cdot,\mu)\|_{\ell^\infty(\Kappa_{m,d})}$}
        \STATE {Compute $\displaystyle \vec{k}_1:=\underset{\vec{k}\in\Kappa_{m,d}}{\textnormal{argmax}}|g(\vec{k},\mu_1)|$}
	\STATE {Set $\displaystyle q^1(\cdot):=\frac{g(\cdot,\mu_1)}{g(\vec{k}_1,\mu_1)}$}
        \STATE {Set $B^1_{11}:=1$}
        \WHILE {$l<Q_{m,d}$}
		\STATE Compute $\displaystyle \mu_{l+1}:=\underset{\mu\in \mathcal{P}_{\text{sample}}}{\textnormal{argmax}}\|(\delta^l g)(\cdot,\mu)\|_{\ell^\infty(\Kappa_{m,d})}$
                \STATE Compute $\displaystyle \vec{k}_{l+1}:=\underset{\vec{k}\in\Kappa_{m,d}}{\textnormal{argmax}}|(\delta^l g)(\vec{k},\mu_{l+1})|$
                \STATE Set $\displaystyle q^{l+1}(\cdot):=\frac{(\delta^l g)(\cdot,\mu_{l+1})}{(\delta^l g)(\vec{k}_{l+1},\mu_{l+1})}$
                \STATE $\displaystyle B^{l+1}_{ij}:=q_j(\vec{k}_{i})$, $1\leq i,j\leq {l+1}$
                \STATE $l\leftarrow l+1$
	\ENDWHILE
\end{algorithmic}
\end{algorithm}

Finally, the online stage of EIM consists in the approximation~\eqref{eq:EIM0}-\eqref{eq:linear_sys} with $l=N^{\rm EIM}$.
Replacing $\beta^l_{l'}(\mu)$ in Equation~\eqref{eq:EIM0} using Equation~\eqref{eq:linear_sys} yields
\begin{equation}
\label{eq:invform}
I^{N^{\rm EIM}}(g)(\vec{k},\mu)=\sum_{l'=1}^{N^{\rm EIM}}\sum_{l''=1}^{N^{\rm EIM}}\left(B^{N^{\rm EIM}}\right)^{-1}_{l',l''}g(\vec{k}_{l''},\mu)q^{l'}(\vec{k}).
\end{equation}
We notice from Algorithm~\ref{algoEIM} that $\underset{1\leq l\leq N^{\rm EIM}}{\rm Span}\left(q^l(\cdot)\right)=\underset{1\leq l\leq N^{\rm EIM}}{\rm Span}\left(g(\cdot,\mu_l)\right)$, and therefore, there exists a matrix $\Gamma\in\mathbb{R}^{N^{\rm EIM}\times N^{\rm EIM}}$ such that, for all $1\leq l\leq N^{\rm EIM}$,
\begin{equation}
\label{eq:relGamma}
\sum_{l'=1}^{N^{\rm EIM}}\Gamma_{l,l'}q^{l'}(\vec{k})=g(\vec{k},\mu_l).
\end{equation}
Replacing $q^{l'}(\vec{k})$ in Equation~\eqref{eq:invform} using Equation~\eqref{eq:relGamma} yields
\begin{equation}
\label{eq:doublesum}
I^{N^{\rm EIM}}(g)(\vec{k},\mu)=\sum_{l'=1}^{N^{\rm EIM}}\sum_{l=1}^{N^{\rm EIM}}\Delta_{l,l'}g(\vec{k}_{l'},\mu)g(\vec{k},\mu_{l}),
\end{equation}
where $\Delta=\left(\Gamma\left(B^{N^{\rm EIM}}\right)^t\right)^{-1}$. From~\cite[Theorem 1.2]{CASENAVE201623}, the matrix $F_{l,l'}=g(\vec{k}_l,\mu_{l'})$, $1\leq l,l'\leq N^{\rm EIM}$ is invertible, and $\Delta=F^{-T}$.
Now denote 
\begin{equation}
\label{eq:lambda}
\lambda^{N^{\rm EIM}}_l(\mu)=\sum_{l'=1}^{N^{\rm EIM}}\Delta_{l,l'}g(\vec{k}_{l'},\mu)
\end{equation}
to obtain
\begin{equation}
\label{eq:EIM}
I^{N^{\rm EIM}}(g)(\vec{k},\mu)= \sum_{l=1}^{N^{\rm EIM}}\lambda^{N^{\rm EIM}}_{l}(\mu)g(\vec{k},\mu_{l}).
\end{equation}

Replacing $g$ in~\eqref{eq:decompA2s} by $I^{N^{\rm EIM}}(g)$ yields
\begin{small}
\begin{equation}
\label{eq:approxAp}
A_\mu^p\approx\sum_{\vec{k}\in\Kappa_{m,d}} \sum_{l=1}^{N^{\rm EIM}}\lambda_l^{N^{\rm EIM}}(\mu)g(\vec{k},\mu_l){T}_{\vec{k},p}
=\sum_{l=1}^{N^{\rm EIM}}\lambda_l^{N^{\rm EIM}}(\mu)\sum_{\vec{k}\in\Kappa_{m,d}}g(\vec{k},\mu_l){T}_{\vec{k},p}=\sum_{l=1}^{N^{\rm EIM}}\lambda_l^{N^{\rm EIM}}(\mu)A_{\mu_l}^p,\quad\mu\in\mathcal{P},~1\leq p\leq m,
\end{equation}
\end{small}
where the last equality is obtained by recognizing $A_{\mu_l}^p$ in Equation~\eqref{eq:decompA2s} at parameter values $\mu_l$.

We obtain the searched expression~\eqref{eq:online_sol}, with $t=N^{\rm EIM}$, and where $\mu_l$ and $\lambda_l(\mu)$ are constructed in respectively the offline and online stages of an EIM on the $g(\vec{k},\mu)$:
\begin{equation}
\label{eq:power_mat_decomp}
A_\mu^p\approx\sum_{l=1}^{N^{\rm EIM}}\lambda_l^{N^{\rm EIM}}(\mu)A_{\mu_l}^p,\quad\mu\in\mathcal{P},~1\leq p\leq m.
\end{equation}

Notice also that the offline stage of EIM involves a sampling of $\mathcal{P}$, and $Q_{m,d}$ indices. Even though $Q_{m,d}$ is independent of the size $N$ of the matrix $A_\mu$, we are limited to moderate values of $m$ and $d$ in practice.

The expression~\eqref{eq:power_mat_decomp} can be used to readily approximate any quantity expressed as a linear evaluation of power of matrices. Yet, quantities such as the inverse matrix or the logarithm of the determinant can be approximated using algorithms involving successive powers of the considered matrix. Hence, we combine in the next section these techniques with the approximation presented in the present section. Under particular conditions, inverses and logarithm of determinants of matrices can be approximated in a nonintrusive fashion using~\eqref{eq:power_mat_decomp}.

\section{Power algorithms}
\label{sec:PowerAlgo}

\subsection{Inverse operators and solution to linear systems}
\label{sec:linsolve}

We recall hereby a classical fixed point results.
\begin{lemma}
Let $A \text{ and }\Psi \in\mathbb{R}^{\mathcal{N}\times \mathcal{N}}$. We consider the following iterative scheme:
\begin{equation}
\label{eq:iteschemeinvmat}
\left\{
\begin{aligned}
X_{0} &= X^0\in \mathbb{R}^{\mathcal{N}\times \mathcal{N}}\\
X_{k+1} &= \left(I-\Psi^{-1}A\right)X_{k} + \Psi^{-1},
\end{aligned}
\right.
\end{equation}
If $\Psi$ is chosen such that $\|I-\Psi^{-1}A\|_2<1$, then the sequence $X_{k}$ converges towards $A^{-1}$ for any initial guess $X^0$.
\end{lemma}

\subsubsection{Sequence approximating the inverse of parametrized matrices}

For concrete implementation, the $k$-th iteration can be evaluated as a series of powers of $A$ and provides an approximation of $A^{-1}$. For a parameter indexed family of matrices, we combine this approximation technique with results of the previous section.

Consider  now a parameter-dependent family of matrices $A_\mu$, $\mu\in\mathcal{P}$,
verifying the affine decomposition~\eqref{eq:decompA0}, and such that, for all $\mu\in\mathcal{P}$, $A_\mu$ is invertible.
We construct a family of approximations of the inverses, $\left(X_{k,\mu}\right)_{k\in \mathbb{N}}$. 
To obtain a uniform convergence with respect to $\mu$, we are led to choose a preconditioner uniform in $\mu$ and a common initial condition: denote
\begin{equation}
\label{eq:bounds}
\left\{
\begin{aligned}
\Psi_0&=\underset{M\in\mathbb{R}^{\mathcal{N}\times \mathcal{N}},~M\rm{ invertible}}{\rm argmin}~\underset{\mu\in\mathcal{P}}{\sup}~\| I-M^{-1}A_\mu\|_2,\\
X_0&=\underset{X\in\mathbb{R}^{\mathcal{N}\times \mathcal{N}}}{\rm argmin}~\underset{\mu\in\mathcal{P}}{\sup}~\|X-A_\mu^{-1}\|_2,
\end{aligned}
\right.
\end{equation}
and let 
\begin{equation}
\label{eq:rho}
\rho=\underset{\mu\in\mathcal{P}}{\sup}\| I-\Psi_0^{-1}A_\mu\|_2,
\end{equation}
and $\epsilon_0=\underset{\mu\in\mathcal{P}}{\sup}\|X_0-A_\mu^{-1}\|_2$. 

We suppose that $\rho<1$ and consider the following iterative scheme:
\begin{equation}
\label{eq:mat_it_scheme}
\left\{
\begin{aligned}
X_{0,\mu} &= X_0\\
X_{k+1,\mu} &= \left(I-\Psi_0^{-1}A_{\mu}\right) X_{k, \mu} + \Psi_0^{-1}.
\end{aligned}
\right.
\end{equation}
An induction shows that:
\begin{equation}
\label{eq:prop_sol}
X_{m,\mu} = \left(I-\Psi_0^{-1}A_{\mu}\right)^m\left(X_0-A_\mu^{-1}\right)+A_\mu^{-1}.
\end{equation}
Taking the norm, we get a uniform bound with respect to $\mu$:
\begin{equation}
\label{eq:rate_conv}
\|X_{m,\mu}-A^{-1}_\mu\|_2\leq \epsilon_0\rho^m,
\end{equation}
which ensures convergence with respect to $m$ since $\rho<1$.

\subsubsection{Powers of a parametrized matrix}
\label{sec:powmatapp}

Define $\alpha_0(\mu)=1$ and $A_0=-\Psi_0$. There holds:
\begin{equation}
\label{newaffinedcomposition}
 \left(I-\Psi_0^{-1}A_{\mu}\right) = -\sum_{l=0}^d \alpha_l(\mu)\Psi_0^{-1}A_l.
\end{equation}
Apply Lemma~\ref{lemmaconstructif} to $\left(I-\Psi_0^{-1}A_{\mu}\right)$ to get
\begin{equation}
\label{eq:powerinv}
 \left(I-\Psi_0^{-1}A_{\mu}\right)^p = \sum_{\vec{k}\in\Kappa_{m,d+1} } \hat{g}\left(\vec{k},\mu\right)\hat{T}_{\vec{k},p},
\end{equation}
where $\hat{g}(\vec{k},\mu)=\prod_{l=0}^d\alpha_l^{k_l}(\mu)$, and $\hat{T}_{\vec{k},p}$ are independent of $\mu$ (the $-$ signs being integrated to the $\hat{T}_{\vec{k},p}$). There holds $\forall\mu\in\mathcal{P}$, $\forall\vec{k}\in\Kappa_{m,d+1}$, $\hat{g}(\vec{k},\mu)=g(c(\vec{k}),\mu)$, where $c$ 
cuts the $0$-th element of $\vec{k}\in \Kappa_{m,d+1}$.
Notice that $c(\vec{k})$ belongs to $\Kappa_{m,d}$. Hence, 
we can apply the EIM approximation to the function~$g(c(\vec{k}),\mu)$ on $\Kappa_{m,d}\times\mathcal{P}$ (see \eqref{eq:EIM}) to obtain :
\begin{equation}
\left(I-\Psi_0^{-1}A_{\mu}\right)^p \approx \sum_{\vec{k}\in\Kappa_{m,d+1}}\left(
\sum_{l=1}^{N^{\rm EIM}}{\lambda_l(\mu) g\left(c(\vec{k}),\mu_l\right)}\right) \hat{T}_{\vec{k},p}.
\end{equation}
We now switch the summations to obtain the desired result:
\begin{equation}
\label{eq:decomp_T}
 \left(I-\Psi_0^{-1}A_{\mu}\right)^p \approx \sum_{l=1}^{N^{\rm EIM}} \lambda_l(\mu)\left(I-\Psi_0^{-1}A_{\mu_l}\right)^p, \quad \mu\in\mathcal{P},\quad 1\leq~p\leq m,
\end{equation}
where we recall that $\mu_l$ and $\lambda_l(\mu)$ are given by EIM on $g(\vec{k},\mu)$. 

Notice that~\eqref{eq:decomp_T} means that we used the same EIM on $g(\vec{k},\mu)$ for the affine approximation~\eqref{eq:decompA0} on $A_\mu$ to get an approximation on $\left(I-\Psi_0^{-1}A_{\mu}\right)^p$.

\subsubsection{Approximation of the inverse of parametrized matrices}

We now go back to the scheme~\eqref{eq:mat_it_scheme} and to show how the approximations \eqref{eq:decomp_T} can be used to approximate the $m$-th approximation $X_{m,\mu}$. Expression~\eqref{eq:prop_sol} is not convenient for this purpose since $A_{\mu}^{-1}$ appears, and we are looking for an efficient approximation of $A_{\mu}^{-1}$. It turns out to be more convenient to consider the following expression obtained by induction:
\begin{equation}
\label{eq:iterative_scheme_2}
X_{m,\mu}=\left(I-\Psi_0^{-1}A_{\mu}\right)^m X_0 + \left(\sum_{k=0}^{m-1} \left(I-\Psi_0^{-1}A_{\mu}\right)^k\right) \Psi_0^{-1},\quad\forall m\in\mathbb{N}.
\end{equation}
To obtain~\eqref{eq:decomp_T}, we use the fact that $\left(I-\Psi_0^{-1}A_{\mu}\right)^p$ depends linearly on $g$, as it explicitly appears in~\eqref{eq:powerinv}. $X_{m,\mu}$ inherits from this linear dependence on $g$; we will make it explicit be denoting now $X_{m,\mu}$ as $X_{m}g_{\mu}$, where $g_{\mu}(\vec{k}):=g(c(\vec{k}),\mu)$.

Now replace the powers of $\left(I-\Psi_0^{-1}A_{\mu}\right)$ in~\eqref{eq:iterative_scheme_2} using~\eqref{eq:decomp_T}:
\begin{equation}
\label{eq:iter_solution}
\begin{aligned}
X_{m,\mu}=X_{m}g_{\mu}&\approx \sum_{l=1}^{N^{\rm EIM}} \lambda_l(\mu)\left(I-\Psi_0^{-1}A_{\mu_l}\right)^m X_0+\sum_{k=0}^{m-1}\sum_{l=1}^{N^{\rm EIM}} \lambda_l(\mu)\left(I-\Psi_0^{-1}A_{\mu_l}\right)^k \Psi_0^{-1}\\
&=\sum_{l=1}^{N^{\rm EIM}} \lambda_l(\mu)\left(\left(I-\Psi_0^{-1}A_{\mu_l}\right)^m X_0+\sum_{k=0}^{m-1}\left(I-\Psi_0^{-1}A_{\mu_l}\right)^k \Psi_0^{-1}\right)\\
&=\sum_{l=1}^{N^{\rm EIM}} \lambda_l(\mu)X_{m}g_{\mu_l}.
\end{aligned}
\end{equation}

The convergence of $X_{m}g_{\mu}$ to $A_{\mu}^{-1}$ with respect to $m$, namely~\eqref{eq:rate_conv}, suggests replacing $X_{m}g_{\mu_l}$ by the inverses $A_{\mu_l}^{-1}$ in~\eqref{eq:iter_solution} and defining
\begin{equation}
\label{eq:iter_solution_2}
\mathcal{X}_\mu^{N^{\rm EIM}}:=\sum_{l=1}^{N^{\rm EIM}}\lambda_l(\mu)A^{-1}_{\mu_l}.
\end{equation}
where $\mathcal{X}_\mu^{N^{\rm EIM}}$ is the obtained approximation of $A^{-1}_\mu$.

\begin{remark}[Nonintrusivity]

In Equation~\eqref{eq:iter_solution_2}, we recall that the coefficients~$\lambda_l(\mu)$ are obtained from an EIM on $g(c(\vec{k}),\mu)$, which only depends on the parametric dependance of~$A_{\mu_l}$, see Equations~\eqref{eq:decompA0} and~\eqref{eq:formula_g}. Therefore, the obtained approximation is nonintrusive in the sense that we only resort to the computation of the quantity of interest (the inverses $A^{-1}_{\mu_l}$) and using some knowledge on the particular form of the problem (the $\alpha_l(\mu)$). In particular, we need to compute neither the $\left(I-\Psi_0^{-1}A_{\mu}\right)^p$ nor the $\hat{T}_{\vec{k},p}$. Notice that the matrices $A_l$ in~\eqref{eq:decompA0}, $1\leq l\leq d$, do not need to be computed either. Moreover, we can compute $A^{-1}_{\mu_l}$ by the method of our choice. Even if the described iterative scheme converges, we can use direct methods
to compute the $A^{-1}_{\mu_l}$, and apply~\eqref{eq:iter_solution_2} to retrieve an approximation of $A^{-1}_\mu$.
We can also compute $A^{-1}_{\mu_l}$ using initial $\mu$-dependent initial guesses $X_{0,\mu}$ and preconditioner $\Psi_{0,\mu}$ in the described iterative scheme.
In particular, we never need to construct $\Psi_0$ and $X_0$, we just need the existence of a matrix $\hat{\Psi}$ such that $\underset{\mu\in\mathcal{P}}{\sup}\|I-\hat{\Psi}A_\mu\|_2<1$.
\end{remark}

\begin{remark}[Solution of linear systems]
Let $b\in\mathbb{R}^{\mathcal{N}}$.
From~\eqref{eq:rate_conv}, there holds $\|X_{k,\mu}b-A^{-1}_\mu b\|_2\leq \rho^k \epsilon_0 \|b\|_2$, which suggests to approximating $A^{-1}_\mu b$ by
\begin{equation}
\label{eq:intersol}
\sum_{l=1}^{N^{\rm EIM}}\lambda_l(\mu)\left(A^{-1}_{\mu_l}b\right).
\end{equation}
\end{remark}

Notice that a key element of the section is the linearity of the function $g\mapsto X_{m}g_{\mu}$, where, from~\eqref{eq:powerinv} and~\eqref{eq:iterative_scheme_2},
\begin{equation}
\label{eq:linearform}
X_{m}g_{\mu}=\sum_{\vec{k}\in\Kappa_{m,d+1}}g(c(\vec{k}),\mu)\left(\hat{T}_{\vec{k},m}X_0+\sum_{l=0}^{m-1}\hat{T}_{\vec{k},l}\Psi_0^{-1}\right).
\end{equation}

\subsection{Logarithm of the determinant}

The logarithm of the determinant (log-det) of a symetric positive definite (SPD) matrix is a quantity receiving interest in the literature. For instance, finding the maximum likelihood estimator of the mean and the covariance matrix of a normal multivariate distribution involves the computation of the log-det of a SPD matrix, see~\cite[Equation (2.2)]{likelihood}.

\subsubsection{Sequence approximating the logarithm of the determinant of parametrized matrices}

Consider a family of parametrized SPD matrices $A_\mu\in\mathbb{R}^{\mathcal{N}\times \mathcal{N}}$ and denote $\rho(A_\mu)$ the spectral radius of $A_{\mu}$. Suppose that $\underset{\mu\in\mathcal{P}}{\sup}~\rho(A_\mu)<\infty$ and that we can determine some $\rho_M>\underset{\mu\in\mathcal{P}}{\sup}~\underset{1\leq i\leq \mathcal{N}}{\max}r_i(A_\mu)$, where $\{r_i(A_\mu)\}_{1\leq i\leq \mathcal{N}}$ denotes the set of eigenvalues of $A_\mu$.
Denote $\rho_0:=\underset{\mu\in\mathcal{P}}{\inf}~\underset{1\leq i\leq \mathcal{N}}{\min}r_i(A_\mu)$, that we suppose strictly positive.
From~\cite[Lemma 5]{logdet},
\begin{equation}
\label{eq:formulalogdet}
\log(\det(A_\mu))=\mathcal{N}\log(\rho_M)-\sum_{k=1}^\infty\frac{{\rm tr}\left((I-\frac{1}{\rho_M}A_\mu)^k\right)}{k}.
\end{equation}
Let $m\in\mathbb{N}$ and consider the following approximation of $\log(\det(A_\mu))-\mathcal{N}\log(\rho_M)$:
\begin{equation}
\label{eq:formulalogdetapprox}
X_{m}g_{\mu}:=-\sum_{k=1}^{m-1}\frac{{\rm tr}\left((I-\frac{1}{\rho_M}A_\mu)^k\right)}{k},
\end{equation}
where we already make explicit the linear dependence in $g$ (see~\eqref{eq:powerinv}).

Since $A_\mu$ is SPD, there exists a family of unitary matrices $U_\mu$ such that $A_\mu=U_\mu D_\mu U_\mu^T$, $D_\mu$ being a diagonal matrix such that ${D_\mu}_{i,i}=r_i(A_\mu)$. From $\left(I-\frac{1}{\rho_M}A_\mu\right)=U_\mu\left(I-\frac{1}{\rho_M}D_\mu\right)U_\mu^T$, there holds ${\rm tr}\left((I-\frac{1}{\rho_M}A_\mu)^k\right)=\sum_{i=1}^\mathcal{N}\left(1-\frac{r_i(A_\mu)}{\rho_M}\right)^k$, from which we infer 
\begin{equation}
\label{eq:ratelogdet}
\begin{aligned}
\left|X_{m}g_{\mu}-\left(\log(\det(A_\mu))-\mathcal{N}\log(\rho_M)\right)\right|&=\sum_{k=m}^{\infty}\frac{{\rm tr}\left((I-\frac{1}{\rho_M}A_\mu)^k\right)}{k}\\
&\leq \mathcal{N}\sum_{k=m}^{\infty}\frac{\left(1-\frac{\rho_0}{\rho_M}\right)^k}{k}\\
&\leq \frac{\mathcal{N}}{m}\sum_{k=m}^{\infty}\left(1-\frac{\rho_0}{\rho_M}\right)^k.
\end{aligned}
\end{equation}
Notice that
\begin{equation}
\label{eq:telescopic_sum}
\begin{aligned}
\frac{\rho_0}{\rho_M}\sum_{k=m}^{\infty}\left(1-\frac{\rho_0}{\rho_M}\right)^k&=\left[1-\left(1-\frac{\rho_0}{\rho_M}\right)\right]\sum_{k=m}^{\infty}\left(1-\frac{\rho_0}{\rho_M}\right)^k\\
&=\sum_{k=m}^{\infty}\left(1-\frac{\rho_0}{\rho_M}\right)^k-\sum_{k=m+1}^{\infty}\left(1-\frac{\rho_0}{\rho_M}\right)^k\\
&=\left(1-\frac{\rho_0}{\rho_M}\right)^m.
\end{aligned}
\end{equation}
Injecting~\eqref{eq:telescopic_sum} in the last term of~\eqref{eq:ratelogdet}, we obtain
\begin{equation}
\label{eq:logdetconv}
\left|X_{m}g_{\mu}+\mathcal{N}\log(\rho_M)-\log(\det(A_\mu))\right|\leq \mathcal{N}\frac{\rho_M}{\rho_0}\frac{\left(1-\frac{\rho_0}{\rho_M}\right)^m}{m},
\end{equation}
which ensures convergence with respect to $m$ since $0\leq 1-\frac{\rho_0}{\rho_M}<1$.

\subsubsection{Powers of a parametrized matrix}

Define $\alpha_0(\mu)=1$ and $A_0=-aI$. There holds:
\begin{equation}
\label{newaffinedcomposition}
\left(I-\frac{1}{a}A_\mu\right) = -\sum_{l=0}^d \frac{\alpha_l(\mu)}{a}A_l.
\end{equation}
We carry out the same analysis as in Section~\ref{sec:powmatapp} to obtain
\begin{equation}
\label{eq:approxpowerB}
\left(I-\frac{1}{a}A_\mu\right)^p\approx\sum_{l=1}^{N^{\rm EIM}}\lambda_l(\mu)\left(I-\frac{1}{a}A_{\mu_l}\right)^p,\quad\mu\in\mathcal{P},\quad 1\leq p\leq m,
\end{equation}
where we recall that $\mu_l$ and $\lambda_l(\mu)$ are given by the EIM on $g(\vec{k},\mu)$. 

\subsubsection{Approximation of the logarithm of the determinant of parametrized matrices}

Replace $\left(I-\frac{1}{a}A_\mu\right)^p$ in the formula~\eqref{eq:formulalogdetapprox} by the right-hand side of~\eqref{eq:approxpowerB} to obtain
\begin{equation}
\label{eq:LDmformula}
\begin{aligned}
X_{m}g_{\mu}&\approx -\sum_{p=1}^{m-1}\frac{1}{p}{{\rm tr}\left(\sum_{l=1}^{N^{\rm EIM}}\lambda_l(\mu)\left(I-\frac{1}{a}A_{\mu_l}\right)^p\right)}\\
&= -\sum_{l=1}^{N^{\rm EIM}}\lambda_l(\mu)\sum_{p=1}^{m-1} \frac{1}{p}{{\rm tr}\left(\left(I-\frac{1}{a}A_{\mu_l}\right)^p\right)}\\
&=\sum_{l=1}^{N^{\rm EIM}}\lambda_l(\mu)X_{m}g_{\mu_l},~\mu\in\mathcal{P}.
\end{aligned}
\end{equation}

Consider the following interpolation property:
\begin{property}[see~Lemma 1 of \cite{Maday}]
\label{interpEIM}
$\forall~1\leq l\leq N^{\rm EIM}$, $\forall\mu\in\mathcal{P}$,
\begin{equation}
I^{N^{\rm EIM}}(g)(\vec{k}_l,\mu)=g(\vec{k}_l,\mu).
\end{equation}
\end{property}
Notice that $\{\vec{k}\in \textnormal{\textlbrackdbl}0;m\textnormal{\textrbrackdbl}^d\textnormal{ such that }|\vec{k}|=0\}=\{\vec{k_0}\}$, where $\vec{k_0}:=(0,0,\cdots 0)$. Besides, $g(\vec{k_0},\mu)=1$ for all $\mu\in\mathcal{P}$. We impose the multi-indice $\vec{k_0}$ to be selected by the EIM in the offline stage, hence the EIM approximation of $g(\vec{k_0},\mu)$ is exact for all $\mu\in\mathcal{P}$ by application of the interpolation Property~\ref{interpEIM}. Hence, $\forall \mu\in\mathcal{P}$, $\sum_{l=1}^{N^{\rm EIM}}\lambda_l(\mu)g(\vec{k_0},\mu_l)=\sum_{l=1}^{N^{\rm EIM}}\lambda_l(\mu)=g(\vec{k_0},\mu)=1$, which enables us to write~\eqref{eq:LDmformula} as
\begin{equation}
\label{eq:LDmformula2}
X_{m}g_{\mu}+\mathcal{N}\log(\rho_M)=\sum_{l=1}^{N^{\rm EIM}}\lambda_l(\mu)\left(X_{m}g_{\mu_l}+\mathcal{N}\log(\rho_M)\right),~\mu\in\mathcal{P}.
\end{equation}

The convergence of $X_{m}g_{\mu}+\mathcal{N}\log(\rho_M)$ to $\log(\det(A_\mu))$ with respect to $m$, namely~\eqref{eq:logdetconv}, suggests replacing $X_{m}g_{\mu_l}+\mathcal{N}\log(\rho_M)$ by $\log(\det(A_{\mu_l}))$ in~\eqref{eq:LDmformula2} and defining
\begin{equation}
\label{eq:iter_solution_3}
\mathcal{X}_\mu^{N^{\rm EIM}}:=\sum_{l=1}^{N^{\rm EIM}}\lambda_l(\mu)\log(\det(A_{\mu_l})),
\end{equation}
where $\mathcal{X}_\mu^{N^{\rm EIM}}$ is the obtained approximation of $\log(\det(A_\mu))$. Notice that we no longer need to compute $\rho_M$, and that any algorithm available to compute $\log(\det(A_{\mu_l}))$, $1\leq l\leq N^{\rm EIM}$, even $\mu$-dependent ones, can be used.

Notice that a key element of the section is the linearity of the function $g\mapsto X_{m}g_{\mu}$, where, from~\eqref{eq:powerinv} and~\eqref{eq:formulalogdetapprox},
\begin{equation}
\label{eq:linearlogdet}
X_{m}g_{\mu}=-\sum_{\vec{k}\in\Kappa_{m,d+1}}g(c(\vec{k}),\mu)\left( \sum_{l=1}^m\frac{{\rm tr}\hat{T}_{\vec{k},l}}{l}\right).
\end{equation}

\subsection{Performance of the approximations}

\subsubsection{Reducibility}
\label{sec:reducibility}

In an industrial context with large-scale computations and a constrained budget, the $N^{\rm EIM}$ in~\eqref{eq:iter_solution_2}, \eqref{eq:intersol}, and \eqref{eq:iter_solution_3} cannot be as large as we want. The success of any nonintrusive procedure will be assessed by the quality of the approximation within the given computation budget. If the approximation yields too large errors, the problem will be considered as nonreducible with the given procedure and the allocated computational budget.
The proposed approximations have been motivated by the iterative schemes~\eqref{eq:iteschemeinvmat} and~\eqref{eq:formulalogdetapprox}, which we recall are not required to be computed in practice. In~\eqref{eq:bounds}, $\Psi_0^{-1}$ can be seen as the best preconditioner uniformly on the parameter space. The problem can be efficiently reduced if this preconditioner is good in the sense that
$\underset{\mu\in\mathcal{P}}{\sup}~\| I-\Psi_0^{-1}A_\mu\|_2=\rho \ll 1$, as can be seen in~\eqref{eq:rate_conv}. 
In high parameter dimension cases, the existence of a good preconditioner is unlikely due to the curse of dimensionality, especially if the interval of variation of each parameter is large. We recall that we do not need to compute $\Psi_0^{-1}$ and just need its existence. The success of the approximation will be assessed \textit{a posteriori}, if a hidden low-rank structure exists, in the same fashion as other \textit{a posteriori} reduced order methods, for instance in the snapshot POD if the eigenvalues of the correlation matrix decrease fast enough. In this context, at given computational budget, we compare our algorithm to some other nonintrusive procedures by computing the approximation errors with respect to reference values in Section~\ref{sec:numerics}.

\subsubsection{Offline cost}

Consider the approximations formulae~\eqref{eq:iter_solution_2}, \eqref{eq:intersol}, and \eqref{eq:iter_solution_3}, that consist in the interpolation of respectively $A_{\mu_l}^{-1}$, $A_{\mu_l}^{-1}b$, and $\log\left(\det\left(A_{\mu_l}\right)\right)$, $1\leq l\leq N^{\rm EIM}$. The construction of these objects is inherent to any nonintrusive approximation method, where the high-fidelity model has to be solved a certain number of times to gather information to derive the approximation. The cost of computing $A_{\mu_l}^{-1}$, $A_{\mu_l}^{-1}b$ and $\log\left(\det\left(A_{\mu_l}\right)\right)$, $1\leq l\leq N^{\rm EIM}$, is then present in any nonintrusive method, and is not related to the offline part of the algorithm derived in the present work. The analysis boils down to assessing the cost of the computation of the coefficients $\lambda_l(\mu)$, $1\leq l\leq N^{\rm EIM}$, in~\eqref{eq:iter_solution_2}, \eqref{eq:intersol}, and \eqref{eq:iter_solution_3}. In our numerical applications, with $d$ imposed by the form of the problem, we determine $m_0$ as the largest $m$ such that $Q_{m,d}=\#\Kappa_{m,d}$ is lower than the computational budget. Then, the offline cost corresponds to the EIM applied to the function $g$ on the sampled spaces $\Kappa_{m_0,d}\times\mathcal{P}_{\rm sample}$. In practice, since the computational budget is constrained (the largest value considered in our numerical experiments for $Q_{m_0,d}$ is 680), we have the opportunity to take a larger sampling of $\mathcal{P}$, which is desired anyway due to the possibly large dimension of $\mathcal{P}$. If the EIM is carried-out until all the $Q_{m_0,d}$ multi-indices in $\Kappa_{m_0,d}$ are selected, the algorithmic complexity is proportional to $Q_{m_0,d}^3\times\#{\mathcal{P}_{\rm sample}}$ : recall that in this case, $Q_{m_0,d}$ corresponds to the number of evaluations of the quantity of interest $A_\mu^{-1}$ or $\log(\det(A_\mu))$. In our numerical experiments, the offline stage of EIM with $\#{\mathcal{P}_{\rm sample}}=10^6$ takes approximately the same time as the construction of the Design Of Experiment (DOE) using MaxProj when comparing with statistical methods, see Section~\ref{sec:numsollinsys} for more details. For instance, with $Q_{m_0,d}=286$ ($m_0=2$ and $d=10$), and $\#{\mathcal{P}_{\rm sample}}=10^6$, both the construction of the DOE and the EIM take approximately 15 minutes.

Notice that in the classical use of EIM for order reduction of general nonlinear models where we want to approach the solution and/or operator, we need to evaluate the high-fidelity model $\#{\mathcal{P}_{\rm sample}}$ number of times: hence a large ${\mathcal{P}_{\rm sample}}$ is not a possile option. However, in the present work, the function $g$ to approximate is known on the complete set $\Kappa_{m_0,d}\times\mathcal{P}$ without solving the high-fidelity model, enabling the possibility of a large ${\mathcal{P}_{\rm sample}}$.

\section{Convergence of the approximation}
\label{sec:convergence}

This section is organized as follows: Section~\ref{settingandnotations} details the setting and notations, Section~\ref{mainresults} states the main results, Section~\ref{technicalproofs} gives the technical proofs, and comments are given in Section~\ref{comments}.

\subsection{Setting}
\label{settingandnotations}

Recall the context of this work: we consider a parameter space $\mathcal{P}$, which is a compact subset of $\mathbb{R}^r$, and denote its Lebesgue mesure by $|\mathcal{P}|$. We also consider a family of matrices  $\{A_\mu\}_{\mu\in\mathcal{P}} \subset \mathbb{R}^{\mathcal{N}\times \mathcal{N}}$. We look for approximations of quantities that can be obtained as limits of power algorithms applied to the matrices $A_\mu$, denoted $\mathcal{L}_\mu$ (standing for "limit" for ease of reading): in the previous section, we considered the inverse matrix: $\mathcal{L}_\mu=A_\mu^{-1}$ and the log-det: $\mathcal{L}_\mu=\log(\det(A_\mu))$.

Let $K$ be a bounded neighborhood of 
$\Kappa_{m,d}$ in $\mathbb{R}^{d}$,
and denote $\mathcal{U}:=L^2(K)$. We recall that $\Kappa_{m,d}=\left\{\vec{k}\in\textnormal{\textlbrackdbl}0;m\textnormal{\textrbrackdbl}^d\textnormal{ such that }|\vec{k}|\leq m\right\}$,
and that $Q_{m,d} = \#\Kappa_{m,d}\leq P_d(m)$, where $P_d(m)$ is a polynomial of degree $d$ in $m$.
Let $\mu\in\mathcal{P}$ and denote $g_{\mu}:\Kappa_{m,d}\ni\vec{k}\mapsto g_{\mu}(\vec{k}):=g(\vec{k},\mu)\in\mathbb{R}$.
Consider the extension of $g_{\mu}$ from $\Kappa_{m,d}$ to $K$: $K\ni\vec{k}\mapsto g_\mu(\vec{k})=\prod_{l=1}^d\alpha_l^{k_l}(\mu)\in\mathbb{R}$, from which we infer $g_{\mu}\in \mathcal{U}=L^2(K)$ -- this point will be important later. We also suppose that $\alpha_l$, $1\leq l\leq d$, are continuous, which ensures the continuity of the functions $\mu\mapsto g(\vec{k},\mu)$ for all $\vec{k}\in K$.

We dispose of a sequence of linear applications $(X_{m})_{m\in\mathbb{N}}\in L\left(\mathcal{U},V\right)$, where $V$ is a Hilbert space of finite dimension $s$ endowed with the scalar product $\left(.,.\right)_{V}$ and its associated norm $\left\|.\right\|_{V}:=\sqrt{\left(.,.\right)_{V}}$ and where $L\left(\mathcal{U},V\right)$ denotes the space of linear applications from $\mathcal{U}$ to $V$. We suppose that the sequence $(X_{m} g_\mu)_{m\in\mathbb{N}}$ converges to $\mathcal{L}_\mu$ in the following sense: for all integer $m$ and all $\mu\in\mathcal{P}$, $\|X_{m} g_\mu-\mathcal{L}_\mu\|_V\leq C_1(m)\underset{m\to\infty}{\longrightarrow}0$. 
We precise here that even if $g_\mu\mapsto X_{m} g_\mu$ is linear, the dependence of the limit $\mathcal{L}_\mu$ with respect to $g_\mu$ is not necessarily linear.
For the inverse matrix and the log-det applications, $\mu\mapsto\mathcal{L}_\mu$ is continuous due to the continuity of the $\alpha_l$, and since $\mathcal{P}$ is a compact subset, $\underset{\mu\in\mathcal{P}}{\sup}\left\|\mathcal{L}_\mu\right\|_V$ can be defined.
Notice that since $\mu\mapsto  X_{m} g_\mu$ is continuous for all $m$, the continuity of $\mu\mapsto \mathcal{L}_\mu$ can be obtained in the general case by assuming the uniform convergence of  $\left(\mu\mapsto  X_{m} g_\mu\right)$ to $\left(\mu\mapsto \mathcal{L}_\mu\right)$ with respect to $m$. We denote $\mathcal{C}^0(\mathcal{P},V)$, the Banach space of the continuous functions from $\mathcal{P}$ to $V$, endowed with the norm $\|w\|_{\mathcal{C}^0(\mathcal{P},V)}=\underset{\mu\in\mathcal{P}}{\sup}{\|w(\mu)\|_{V}}$.

For the inverse operators, the linear operator in $g_\mu$ is
$$\mathcal{U}\ni g_\mu\mapsto X_{m} g_\mu:=\sum_{\vec{k}\in\Kappa_{m,d+1}}g(c(\vec{k}),\mu)\left(\hat{T}_{\vec{k},m}X_0+\sum_{l=0}^{m-1}\hat{T}_{\vec{k},l}\Psi_0^{-1}\right)\in V:=\mathbb{R}^{\mathcal{N}\times\mathcal{N}},$$
hence $s=\mathcal{N}^2$, the limit is $\mathcal{L}_\mu:=A_\mu^{-1}$, and $C_{1}(m)=\epsilon_0\rho^{m}$, see~\eqref{eq:rate_conv}.

For the logarithm of the determinant, the linear operator in $g_\mu$ is
$$\mathcal{U}\ni g_\mu\mapsto X_{m} g_\mu:=-\sum_{\vec{k}\in\Kappa_{m,d+1}}g(c(\vec{k}),\mu)\left( \sum_{l=1}^m\frac{{\rm tr}\hat{T}_{\vec{k},l}}{l}\right)\in V:=\mathbb{R},$$
hence $s=1$,  the limit is $\mathcal{L}_\mu:=\log(\det(A_\mu))-\mathcal{N}\log(\rho_M)$, and $C_{1}(m)=\mathcal{N}\frac{\rho_M}{\rho_0}\frac{\left(1-\frac{\rho_0}{\rho_M}\right)^{m}}{m}$, see~\eqref{eq:logdetconv}.

Consider the following EIM approximation of $g$: 
$$I^{N^{\rm EIM}}(g)(\vec{k},\mu):=\sum_{l=1}^{N^{\rm EIM}}\lambda_l(\mu)g(\vec{k},\mu_l),\quad \vec{k}\in\Kappa_{m,d},\quad\mu\in\mathcal{P},$$ where we recall that
$\lambda_l(\mu)=\sum_{l'=1}^{N^{\rm EIM}}\Delta_{l,l'}g(\vec{k}_{l'},\mu)$, $\Delta=(F)^{-T}$ where $F_{l,l'}=g(\vec{k}_l,\mu_{l'})$, $1\leq l,l'\leq N^{\rm EIM}$, where $\vec{k}_l$ and $\mu_{l'}$ are selected during the offline stage of EIM. We denote 
\begin{equation}
\label{eq:defdeltan}
\delta_{N^{\rm EIM}}=\left\|I^{N^{\rm EIM}}(g)-g\right\|_{L^2(\mathcal{P},\mathcal{U})}:=\sqrt{\int_{\mu\in\mathcal{P}}\left\|I^{N^{\rm EIM}}(g)(\cdot,\mu)-g(\cdot,\mu)\right\|^2_{\mathcal{U}}},
\end{equation}
which can be defined thanks to the compactness of $\mathcal{P}$ and the continuity of $\mu\mapsto g(\vec{k},\mu)$ for all $\vec{k}\in K$, ensuring also the continuity of $\mu\mapsto\lambda_l(\mu)$, $1\leq l\leq N^{\rm EIM}$, yielding the intregrability. 
Denote $\hat{Q}_{m,d}\leq Q_{m,d}$, the rank of the matrix $\left(g(\vec{k}_i,\mu_{j})\right)_{i,j}$, $1\leq i\leq Q_{m,d}$, $1\leq j\leq\#\mathcal{P}$. 
Owing to Property~\ref{interpEIM}, $N^{\rm EIM}=\hat{Q}_{m,d}$ implies that the EIM approximation is exact on the whole domain $\Kappa_{m,d} \times \mathcal{P}$. Hence, we now consider values for $N^{\rm EIM}$ smaller than $\hat{Q}_{m,d}$. For ease of reading, we set $N := N^{\rm EIM}$, keeping in mind the dependency of $N$ in $m$.

In what follows, we denote by $(.,.)_{\mathcal{U}}$ the scalar product on $\mathcal{U}$ and $\|.\|_{\mathcal{U}}$ its associated norm. The corresponding inner product is the $L^{2}-$ one. As explained at the beginning of the section, $g_{\mu}\in \mathcal{U}$, for all $\mu\in\mathcal{P}$.
Denote the set $S=\left\{g_{\mu_l}\right\}_{1\leq l\leq N}\subset\mathcal{U}$ where the $\mu_l$ are the parameter values selected by the EIM on $g$. We apply the POD technique to the set $S$, see Table~\ref{tablePOD} for the obtained properties and~\cite{sirovich,bergmann2004optimisation} for more details and justifications.

\def\arraystretch{1.8}
\begin{table}[H]
\centering
\begin{tabular}{|l|l|}
\hline
 set & $S=\left\{g_{\mu_l}\right\}_{1\leq l\leq N}$ \\ \hline
correlation operator  &  $C_{pq} = \left(g_{\mu_p},g_{\mu_q}\right)_{\mathcal{U}}$  \\ \hline
eigenvalue problem  & $\tau_n\xi_{n,p}=\frac{1}{N}\sum_{q=1}^{N}C_{pq}\xi_{n,q}$  \\ \hline
POD modes & $\Phi_n=\frac{1}{\sqrt{N\tau_n}}\sum_{p=1}^N\xi_{n,p}g_{\mu_p}$  \\ \hline
eigenvalues property &  $\tau_n = \frac{1}{N}\sum_{p=1}^N\left(g_{\mu_p},\Phi_n\right)^2_{\mathcal{U}}$ \\ \hline
eigenfunctions orthonormality  & $\sum_{p=1}^N\xi_{n,p}\xi_{m,p}=\delta_{n,m}$ \\ \hline
POD modes orthonormality &  $\left(\Phi_n,\Phi_m\right)_{\mathcal{U}}=\delta_{n,m}$ \\ \hline
\end{tabular}
\caption{Definitions and properties resulting from the POD on the sets $S$\label{tablePOD}}
\end{table}

The approximation of $\mathcal{L}_\mu$, denoted $\mathcal{X}_\mu^{N}$, is defined as
\begin{equation}\label{eq:defxronde}
\mathcal{X}_\mu^{N}:=\sum_{n=1}^{N}\lambda_n(\mu)\mathcal{L}_{\mu_n}.
\end{equation}

\subsection{Main results}
\label{mainresults}

In this section, we give two different bounds for the error made by the approximation $\mathcal{X}_\mu^{N}$ of  $\mathcal{L}_\mu$ the first one
involves a rather abstract vector space, the second one makes use of the relation between the functions $g_\mu$, on which the EIM approximation is carried out, and the approximated object $\mathcal{L}_\mu$, through the iterative schemes $X_m g_\mu$.

Define $\displaystyle Z_N:=\frac{1}{N}\sum_{n=1}^N\frac{{\int_{\mu\in\mathcal{P}}\left(g_{\mu},\Phi_n\right)^2_{\mathcal{U}}}}{\frac{|\mathcal{P}|}{N}\sum_{p=1}^N\left(g_{\mu_p},\Phi_n\right)^2_{\mathcal{U}}}$:  in the quotient, the denominator is an approximation of the numerator, leading to the boundedness of $(Z_N)_N$.
Define also $\mathcal{S}^{sN}$ the smallest sN-dimensional subspace of $\mathcal{C}^0(\mathcal{P},V)$ containing the image of the application $v\mapsto\mathcal{J}^Nv$, defined by $\forall\mu\in\mathcal{P}$, $(\mathcal{J}^Nv)(\mu):=\sum_{n=1}^N {\left(g_{\mu},\Phi_n\right)_{\mathcal{U}}}\sum_{p=1}^NG^{-1}_{np}v(\mu_p)$, where $G_{np}=\left(g_{\mu_n},\Phi_p\right)_{\mathcal{U}}\in\mathbb{R}^{N\times N}$ is an invertible matrix. The boundedness of $(Z_N)_N$, the dimension of $\mathcal{S}^{sN}$ and the invertibility of $G$ will be justified in Section~\ref{technicalproofs}.
\begin{proposition}\label{convergence1}
For any integer $m$ and $1\leq N<\hat{Q}_{m,d}$,
\begin{equation}\label{eq:prop6}
\frac{1}{|\mathcal{P}|}\int_{\mu\in\mathcal{P}}\|\mathcal{X}_{\mu}^{N}-\mathcal{L}_\mu\|^{2}_V\leq{
4\left(1+N^2{Z_N}\right)\left(\theta^{sN}_{\mathcal{L}}\right)^2+\frac{8}{|\mathcal{P}|}{\underset{\mu\in\mathcal{P}}{\sup}\left\|\mathcal{L}_\mu\right\|_V^2}\frac{N\delta_{N}^2}{\tau_N}},
\end{equation}
where
\begin{equation}\label{eq:theta}
\theta_{\mathcal{L}}^{sN}:=\underset{\varphi\in\mathcal{S}^{sN}}{\inf}~\underset{\mu\in\mathcal{P}}{\sup}~{\left\|\mathcal{L}_{\mu}-\varphi(\mu)\right\|_V}.
\end{equation}

In the case $N=\hat{Q}_{m,d}$ where the EIM approximation is exact: for any integer $m$
\begin{equation}\label{eq:prop6b}
\frac{1}{|\mathcal{P}|}\int_{\mu\in\mathcal{P}}\|\mathcal{X}_{\mu}^{\hat{Q}_{m,d}}-\mathcal{L}_\mu\|^{2}_V\leq 
4\left(1+{\hat{Q}_{m,d}}^2{Z_{\hat{Q}_{m,d}}}\right)\left(\theta^{s{\hat{Q}_{m,d}}}_{\mathcal{L}}\right)^2.
\end{equation}
\end{proposition}

\begin{proposition}\label{convergence3}
For any integer $m$ and $1\leq N<\hat{Q}_{m,d}$,
\begin{equation}\label{eq:prop7}
\frac{1}{|\mathcal{P}|}\int_{\mu\in\mathcal{P}}\|\mathcal{X}_{\mu}^{N}-\mathcal{L}_\mu\|^{2}_V\leq 
3C_1^2(m)\left(1+2N^2{Z_{N}}+8\frac{1}{|\mathcal{P}|}\frac{N\delta_N^2}{\tau_N}\right)+3\frac{1}{|\mathcal{P}|}\int_{\mu\in\mathcal{P}}\left\| X_m\left(I^{N}(g_\mu)-g_\mu\right)\right\|_V^2.
\end{equation}

In the case $N=\hat{Q}_{m,d}$ where the EIM approximation is exact: for any integer $m$
\begin{equation}\label{eq:prop7b}
\frac{1}{|\mathcal{P}|}\int_{\mu\in\mathcal{P}}\|\mathcal{X}_{\mu}^{\hat{Q}_{m,d}}-\mathcal{L}_\mu\|^{2}_V\leq 
3C_1^2(m)\left(1+2\hat{Q}_{m,d}^2{Z_{\hat{Q}_{m,d}}}\right).
\end{equation}
\end{proposition}

\begin{remark}
The bounds in~\eqref{eq:prop6} and~\eqref{eq:prop7} involve $\frac{N\delta_{N}^2}{\tau_N}$ and  $\left\| X_m\left(I^{N}(g_\mu)-g_\mu\right)\right\|_V^2$, which are difficult to describe: on the one hand the asymptotic behavior of $\frac{N\delta_{N}^2}{\tau_N}$ exibits an indeterminate form, and on the other hand the operator norm of $X_m$ is hard to estimate. However, thanks to the interpolation property of the EIM, we know that $\delta_{\hat{Q}_{m,d}}=0$ and $I^{\hat{Q}_{m,d}}(g_\mu)=g_\mu$, while $\tau_{\hat{Q}_{m,d}}>0$: sharper upper bounds are derived in this particular case. 

In~\eqref{eq:prop6b}, $Z_{\hat{Q}_{m,d}}$ is bounded, and $\theta^{s{\hat{Q}_{m,d}}}_{\mathcal{L}}$ is related to a certain Kolmogorov width, which are usually assumed to compensate for exponential or polynomial growth in approximation problems; in our case, we only need to assume the convergence of $\hat{Q}_{m,d}\theta^{s{\hat{Q}_{m,d}}}_{\mathcal{L}}$. This is thoroughly commented in Section~\ref{comments}. In~\eqref{eq:prop7b}, the convergence is ensured by the properties of the considered power algorithm through $C_1(m)$, as we explicit in the following corollary.
\end{remark}

\begin{coro}\label{coro1}
In the case $N = \hat{Q}_{m,d}$ where the EIM approximation is exact, the bound of Proposition~\ref{convergence3} is
\begin{itemize}
\item for the inverse matrix :
\begin{equation}
\label{eq:coro1}
\frac{1}{|\mathcal{P}|}\int_{\mu\in\mathcal{P}}\|\mathcal{X}_{\mu}^{\hat{Q}_{m,d}}-A_\mu^{-1}\|^{2}_V\leq 3\epsilon_0^2\rho^{2m}\left(1+2P^2_{d}(m){Z_{\hat{Q}_{m,d}}}\right),
\end{equation}
\item for the log-det :
\begin{equation}
\frac{1}{|\mathcal{P}|}\int_{\mu\in\mathcal{P}}|\mathcal{X}_{\mu}^{\hat{Q}_{m,d}}-\log(\det(A_\mu))|^{2}\leq3\frac{\mathcal{N}^2}{m^2}\frac{\rho_M^2}{\rho_0^2}{\left(1-\frac{\rho_0}{\rho_M}\right)^{2m}}\left(1+2P^2_{d}(m){Z_{\hat{Q}_{m,d}}}\right),
\end{equation}
\end{itemize}
where we recall that $P_d(m)$ is a polynomial of degree $d$ in $m$, such that $\hat{Q}_{m,d}\leq {Q}_{m,d}\leq P_d(m)$.
The approximation converges under the condition that $\rho<1$ for the inverse matrix, and that $1-\frac{\rho_0}{\rho_M}<1$ for the log-det.
\end{coro}

\begin{remark}[Reducibility]
In the case of the inverse matrix, the bound in~\eqref{eq:coro1} converges under the strong assumption that $\rho<1$, where we recall that $\rho=\underset{\mu\in\mathcal{P}}{\sup}\| I-\Psi_0^{-1}A_\mu\|_2$. The strength of the assumption lies in the existence of a good preconditioner $\Psi_0$ uniformly on the possibly large dimensional parameter space $\mathcal{P}$, which we have related to the reducibility of the problem at hand in Section~\ref{sec:reducibility}.
\end{remark}

\subsection{Technical proofs}
\label{technicalproofs}

In this section, we start by giving two results on the POD basis $\left(\Phi_n\right)_{n\in\mathbb{N}}$: Intermediate Result~\ref{interm2} and~\ref{interm1}, from which we derive the proofs of Proposition~\ref{convergence1} and~\ref{convergence3}.

\begin{interm_result}\label{interm2}
The matrix $G_{np}=\left(g_{\mu_n},\Phi_p\right)_{\mathcal{U}}\in\mathbb{R}^{N\times N}$ is invertible, and $\forall\vec{r}\in\mathbb{R}^N$, $\forall\vec{v}:=(v_n)_{1\leq n\leq N}$ with $v_n\in V$, there holds
\begin{equation}
\left\|\vec{r}~^t\cdot\left(G^{-1}\vec{v}\right)\right\|_{V}^2\leq \left(\sum_{n=1}^N\|v_n\|_V^2\right)\left(\sum_{n=1}^N\frac{r_n^2}{\tau_n}\right).
\end{equation}
\end{interm_result}

\renewcommand*{\proofname}{\textbf{Proof of Intermediate result~\ref{interm2}}}
\begin{proof}
The family $\left( g_{\mu_n}\right)_{1\leq n\leq N}$ is free over $\mathcal{U}$. Indeed, let $\vec{a}\in \mathbb{R}^N$ such that $\displaystyle\sum^{N}_{n=1}a_n g_{\mu_n}=0$. The equality holds in particular for the indices $\vec{k}_l\in\Kappa_{m,d}$ selected by EIM: $\forall~ 1\leq l\leq N$, $\displaystyle\sum^{N}_{n=1}a_n g(\mu_n,\vec{k}_l)=0$. By construction of the EIM, the matrix $(g(\vec{k}_l,\mu_n))_{1\leq l,n\leq N}$ is invertible (see \cite[Lemma~2.2]{CASENAVE201623}), in particular the rows of this matrix form a free family. This entails that all the $a_n$ are zero, which proves that the family $\left( g_{\mu_n}\right)_{1\leq n\leq N}$ is free over $\mathcal{U}$. Then, thanks to the orthonormality of the basis $\left(\Phi_n\right)_{n\in\mathbb{N}}$, $g_{\mu_n}=\displaystyle\sum^{N}_{p=1}\left(g_{\mu_n},\Phi_p\right)_{\mathcal{U}}\Phi_p=\displaystyle\sum^{N}_{p=1}G_{np}\Phi_p$, for all $1\leq  n\leq N$, which means that $G$ is the change of basis matrix from $\left(\Phi_n\right)_{1\leq n\leq N}$ to $\left(g_{\mu_p}\right)_{1\leq p\leq N}$, and is therefore invertible. 

Hence, $\Phi_n=\sum_{p=1}^{N}G_{np}^{-1}g_{\mu_p}$. From the definition of the POD modes $\Phi_n$ (see Table~\ref{tablePOD}) and due to the fact that $\left(g_{\mu_i}\right)_{1\leq i\leq N}$ is a set of linear independent vectors, we obtain $G_{np}^{-1}=\frac{1}{\sqrt{N\tau_n}}\xi_{n,p}$ for all $1\leq n,p\leq N$. Then, using the eigenfunctions orthonormality  (see Table~\ref{tablePOD}), 
\begin{equation}
\label{eq:sumGm1_2}
\sum_{p=1}^N\left(G_{np}^{-1}\right)^2=\frac{1}{N\tau_n}\sum_{p=1}^N\xi_{n,p}^2=\frac{1}{N\tau_n}.
\end{equation}
There holds
\begin{subequations}
\begin{align}
\left\|\vec{r}~^t\cdot\left(G^{-1}\vec{v}\right)\right\|_{V}^2&=\left(\left\|\sum_{n=1}^Nr_n\sum_{p=1}^NG_{np}^{-1}v_p\right\|_V\right)^2\label{eq:r1}\\
&\leq\left(\sum_{n=1}^N|r_n|\left\|\sum_{p=1}^NG_{np}^{-1}v_p\right\|_V\right)^2\label{eq:r2}\\
&\leq N\sum_{n=1}^Nr_n^2\left\|\sum_{p=1}^NG_{np}^{-1}v_p\right\|_V^2\label{eq:r3}\\
&\leq N\sum_{n=1}^Nr_n^2\left(\sum_{p=1}^N\left|G_{np}^{-1}\right|\left\|v_p\right\|_V\right)^2\label{eq:r4}\\
&\leq N\sum_{n=1}^Nr_n^2\left(\sum_{p=1}^N\left(G_{np}^{-1}\right)^2\right)\left(\sum_{p=1}^N\left\|v_p\right\|_V^2\right)\label{eq:r5}\\
&\leq \left(\sum_{n=1}^N\|v_n\|_V^2\right)\left(\sum_{n=1}^N\frac{r_n^2}{\tau_n}\right),\label{eq:r6}
\end{align}
\end{subequations}
where the Jensen inequality is applied to the square function between~\eqref{eq:r2} and~\eqref{eq:r3}, the Cauchy-Schwarz inequality is applied between~\eqref{eq:r4} and~\eqref{eq:r5}, and where~\eqref{eq:r6} is obtained from~\eqref{eq:r5} using~\eqref{eq:sumGm1_2}, which ends the proof.
\end{proof}

\begin{remark}
In the proof of Intermediate result~\ref{interm2}, we could have simply introduced an orthonormal basis $(\Phi_n)_{1\leq n\leq N}$, obtained for instance from a Gram-Schmidt orthonormalization of the family $g_{\mu_l}$ $1\leq l\leq N$, and used $\sum_{p=1}^N\left(G_{np}^{-1}\right)^2\leq \|G^{-1}\|^2_F$, where $\|G^{-1}\|_F$ is the Frobenius norm of $G^{-1}$, depending on $N$. However, doing so would not have yielded a bound with an explicit dependence on $N$.
\end{remark}

\begin{interm_result}\label{interm1}
$$\int_{\mu\in\mathcal{P}}\left\|\Pi^{\Phi}_{N}I^{N}(g_{\mu})-\Pi^{\Phi}_{N}g_{\mu}\right\|^{2}_{\mathcal{U}}\leq 4\delta_N^2,$$
where $\Pi^{\Phi}_{N}$ is the orthogonal projection operator onto the subspace $\underset{1\leq n\leq N}{\rm Span}(\Phi_n)$ and where we recall that $\delta_N$ quantifies the EIM approximation error, see~\eqref{eq:defdeltan}.
\end{interm_result}
\renewcommand*{\proofname}{\textbf{Proof of Intermediate result~\ref{interm1}}}
\begin{proof}
$$\Pi^{\Phi}_{N}I^{N}(g_{\mu})-\Pi^{\Phi}_{N}g_{\mu}= \Pi^{\Phi}_{N}I^{N}(g_{\mu})-g_{\mu}+g_{\mu}-\Pi^{\Phi}_{N}g_{\mu}.$$
Recall that $\underset{1\leq n\leq N}{\rm Span}(\Phi_n)=\underset{1\leq n\leq N}{\rm Span}(g_{\mu_n})$ from Intermediate Result~\ref{interm2}, providing $\Pi^{\Phi}_{N}I^{N}(g_{\mu})=I^{N}(g_{\mu})$.
Using the triangular inequality hence yields 
\begin{eqnarray*}
\int_{\mu\in\mathcal{P}}\left\|\Pi^{\Phi}_{N}I^{N}(g_{\mu})-\Pi^{\Phi}_{N}g_{\mu}\right\|^{2}_{\mathcal{U}}&\leq&2\int_{\mu\in\mathcal{P}}\left\|I^{N}(g_{\mu})-g_{\mu}\right\|^{2}_{\mathcal{U}}+2\int_{\mu\in\mathcal{P}}\left\|g_{\mu}-\Pi^{\Phi}_{N}g_{\mu}\right\|^{2}_{\mathcal{U}}\\
&=&2\int_{\mu\in\mathcal{P}}\left\|I^{N}(g_{\mu})-g_{\mu}\right\|^{2}_{\mathcal{U}}+2\int_{\mu\in\mathcal{P}}\underset{v\in\underset{1\leq n\leq N}{\rm Span}(\Phi_n)}{\inf}\left\|g_{\mu}-v\right\|^{2}_{\mathcal{U}}\\
&\leq&4\int_{\mu\in\mathcal{P}}\left\|I^{N}(g_{\mu})-g_{\mu}\right\|^{2}_{\mathcal{U}}\\
&\leq&4\delta_N^2,
\end{eqnarray*}
which ends the proof.
\end{proof}

\renewcommand*{\proofname}{\textbf{Proof of Proposition~\ref{convergence1}}}
\begin{proof}
First, we recall the definition of the application $\mathcal{C}^0(\mathcal{P},V)\ni v\mapsto\mathcal{J}^Nv\in \mathcal{C}^0(\mathcal{P},V)$ such that $\forall\mu\in\mathcal{P}$, $(\mathcal{J}^Nv)(\mu):=\sum_{n=1}^N {\left(g_{\mu},\Phi_n\right)_{\mathcal{U}}}\sum_{p=1}^NG^{-1}_{np}v(\mu_p)$, and of the subspace $\mathcal{S}^{sN}$: it is the smallest subspace of $ \mathcal{C}^0(\mathcal{P},V)$ containing the image of $\mathcal{J}^N$. To see that the dimension of $\mathcal{S}^{sN}$ is $sN$, denote $(e_i)_{1\leq i\leq s}$ a basis of $V$; any element $w$ of $V$ can be expressed in this basis as follows: $w=\sum_{i=1}^s\eta_i(w)e_i$. Then, for all $v\in\mathcal{C}^0(\mathcal{P},V)$, there holds $(\mathcal{J}^N v)(\mu)=\sum_{p=1}^N\sum_{i=1}^s\left(\eta_i(v(\mu_p))\right)\left(\sum_{n=1}^N {\left(g_{\mu},\Phi_n\right)_{\mathcal{U}}}G^{-1}_{np}e_i\right)$, where $\eta_i(v(\mu_p))\in\mathbb{R}$ and $\sum_{n=1}^N {\left(g_{\mu},\Phi_n\right)_{\mathcal{U}}}G^{-1}_{np}e_i\in \mathcal{C}^0(\mathcal{P},V)$ and is independent of $v$. This proves that the family $\left(\left(\sum_{n=1}^N {\left(g_{\mu},\Phi_n\right)_{\mathcal{U}}}G^{-1}_{np}\right)e_i\right)_{1\leq p\leq N,~1\leq i\leq s}$ is composed of spanning vectors of $\mathcal{S}^{sN}$. We conclude by noticing that this family is free in $\mathcal{C}^0(\mathcal{P},\mathbb{R})\times V\subset  \mathcal{C}^0(\mathcal{P},V)$: $(e_i)_{1\leq i\leq s}$ is a basis of $V$, and let $\omega\in\mathbb{R}^N$ such that $\sum_{p=1}^N\omega_p\left(\sum_{n=1}^N {\left(g_{\mu},\Phi_n\right)_{\mathcal{U}}}G^{-1}_{np}\right)=0$ in $\mathcal{C}^0(\mathcal{P},\mathbb{R})$. In particular, $\forall 1\leq q\leq N$, $\sum_{p=1}^N\omega_p\left(\sum_{n=1}^N {\left(g_{\mu_q},\Phi_n\right)_{\mathcal{U}}}G^{-1}_{np}\right)=\sum_{p=1}^N\omega_p\left(GG^{-1}\right)_{qp}=\omega_q=0$.

We now go back to the control of $$\int_{\mu\in\mathcal{P}}\left\|\mathcal{X}_\mu^N-\mathcal{L}_\mu\right\|_V^2=\int_{\mu\in\mathcal{P}}\left\|\sum_{n=1}^N\lambda_n(\mu)\mathcal{L}_{\mu_n}-\mathcal{L}_\mu\right\|_V^2.$$
Denote $\vec{\lambda}(\mu):=\left(\lambda_n(\mu)\right)_{1\leq n\leq N}$, $\vec{\mathcal{L}}:=\left(\mathcal{L}_{\mu_n}\right)_{1\leq n\leq N}$ and $\displaystyle\vec{h}(\mu)=\left({\left(g_\mu,\Phi_n\right)_{\mathcal{U}}}\right)_{1\leq n\leq N}$. With these notations, $(\mathcal{J}^N\mathcal{L})(\mu):=\vec{h}^t(\mu)\cdot G^{-1}\vec{\mathcal{L}}$. Then,
\begin{subequations}
\begin{align}
\int_{\mu\in\mathcal{P}}\left\|\sum_{n=1}^N\lambda_n(\mu)\mathcal{L}_{\mu_n}-\mathcal{L}_\mu\right\|_V^2&=\int_{\mu\in\mathcal{P}}\left\|\vec{\lambda}^t(\mu)\cdot\vec{\mathcal{L}}-\mathcal{L}_\mu\right\|_V^2\\
&\leq 2\int_{\mu\in\mathcal{P}}\left\|\left(G^t\vec{\lambda}(\mu)-\vec{h}(\mu)\right)^t\cdot\left(G^{-1}\vec{\mathcal{L}}\right)\right\|_V^2+2\int_{\mu\in\mathcal{P}}\left\|(\mathcal{J}^N\mathcal{L})(\mu)-\mathcal{L}_{\mu}\right\|_V^2.\label{eq:propc}
\end{align}
\end{subequations}
We now control the two terms in~\eqref{eq:propc}:
\begin{itemize}
\item(first term in~\eqref{eq:propc})
\begin{subequations}
\begin{align}
2\int_{\mu\in\mathcal{P}}\left\|\left(G^t\vec{\lambda}(\mu)-\vec{h}(\mu)\right)^t\cdot\left(G^{-1}\vec{\mathcal{L}}\right)\right\|_V^2
&\leq 2\left(\sum_{n=1}^N\mathcal{L}_{\mu_n}^2\right)\int_{\mu\in\mathcal{P}}\sum_{n=1}^N\frac{1}{\tau_n}\left(G^t\vec{\lambda}(\mu)-\vec{h}(\mu)\right)^2_n\label{eq:control6a}\\
&\leq 2\underset{\mu\in\mathcal{P}}{\sup}\left\|\mathcal{L}_\mu\right\|_V^2N\int_{\mu\in\mathcal{P}}\sum_{n=1}^N\frac{1}{\tau_n}\left(\sum_{p=1}^N\lambda_p(\mu)g_{\mu_p}-g_{\mu},\Phi_n\right)_{\mathcal{U}}^2\label{eq:control6b}\\
&\leq2\underset{\mu\in\mathcal{P}}{\sup}\left\|\mathcal{L}_\mu\right\|_V^2\frac{N}{\tau_N}\int_{\mu\in\mathcal{P}}\left\|\Pi_N^\Phi\left(I^Ng_\mu-g_\mu\right)\right\|^2_{\mathcal{U}}\label{eq:control6c}\\
&\leq 8\underset{\mu\in\mathcal{P}}{\sup}\left\|\mathcal{L}_\mu\right\|_V^2\frac{N\delta_N^2}{\tau_N},\label{eq:control6d}
\end{align}
\end{subequations}
where we applied Intermediate result~\ref{interm2} to  $\vec{r}=G^t\vec{\lambda}(\mu)-\vec{h}(\mu)$ and $\vec{v}=\vec{\mathcal{L}}$ in~\eqref{eq:control6a} and Intermediate result~\ref{interm1} between~\eqref{eq:control6c} and~\eqref{eq:control6d}.
\item(second term in~\eqref{eq:propc})
For all $v\in\mathcal{C}^0(\mathcal{P},V)$, the application $(\mathcal{J}^Nv)(\mu)$ defines an interpolation, since $\forall~1\leq q\leq N$,
\begin{equation}
\label{eq:interpprop}
\left(\mathcal{J}^Nv\right)(\mu_q)=\sum_{n=1}^N {\left(g_{\mu_q},\Phi_n\right)_{\mathcal{U}}}\sum_{p=1}^NG^{-1}_{np}v(\mu_p)=\sum_{p=1}^N\left(GG^{-1}\right)_{qp}v(\mu_p)=v(\mu_q).
\end{equation}
Moreover, $\mathcal{J}^N$ is a linear projector onto $\mathcal{S}^{sN}$ since using~\eqref{eq:interpprop}, for all $v\in\mathcal{C}^0(\mathcal{P},V)$,
\begin{equation}
\label{eq:projprop}
\left(\mathcal{J}^N\mathcal{J}^Nv\right)(\mu)=\sum_{n=1}^N {\left(g_{\mu},\Phi_n\right)_{\mathcal{U}}}\sum_{p=1}^NG^{-1}_{np}(\mathcal{J}^Nv)(\mu_p)=\sum_{n=1}^N {\left(g_{\mu},\Phi_n\right)_{\mathcal{U}}}\sum_{p=1}^NG^{-1}_{np}v(\mu_p)=\left(\mathcal{J}^Nv\right)(\mu).
\end{equation}
Now denote $\varphi_0=\underset{\varphi\in\mathcal{S}^{sN}}{\rm arginf}~{\|\varphi-\mathcal{L}\|_{\mathcal{C}^0(\mathcal{P},V)}}$. Since $\mathcal{J}^N$ is a projector, $\varphi_0=\mathcal{J}^N\varphi_0$ in $\mathcal{C}^0(\mathcal{P},V)$. We control the second term in~\eqref{eq:propc} with
\begin{equation}
\begin{aligned}
2\int_{\mu\in\mathcal{P}}\left\|(\mathcal{J}^N\mathcal{L})(\mu)-\mathcal{L}_{\mu}\right\|_V^2 =2\left\|\mathcal{J}^N \mathcal{L}-\mathcal{L}\right\|^2_{L^2(\mathcal{P},V)}&\leq4\left\|\mathcal{J}^N \mathcal{L}-\varphi_0\right\|^2_{L^2(\mathcal{P},V)}+4\left\|\varphi_0-\mathcal{L}\right\|^2_{L^2(\mathcal{P},V)}\\
&\leq 4\left\|\mathcal{J}^N\left(\mathcal{L}-\varphi_0\right)\right\|^2_{L^2(\mathcal{P},V)}+4|\mathcal{P}|\left\|\varphi_0-\mathcal{L}\right\|^2_{\mathcal{C}^0(\mathcal{P},V)}\\
&=4\int_{\mu\in\mathcal{P}}\left\|\vec{h}^t(\mu)\cdot G^{-1}\left(\overrightarrow{\mathcal{L}-\varphi_0}\right)\right\|_V^2+4|\mathcal{P}|\left(\theta^{sN}_{\mathcal{L}}\right)^2.
\end{aligned}
\end{equation}
Using Intermediate result~\ref{interm2} with $\vec{r}=\vec{h}(\mu)$, there holds,
\begin{subequations}
\begin{align}
4\int_{\mu\in\mathcal{P}}\left\|\vec{h}^t(\mu)\cdot G^{-1}\left(\overrightarrow{\mathcal{L}-\varphi_0}\right)\right\|_V^2
&\leq 4\int_{\mu\in\mathcal{P}}\left(\sum_{n=1}^N\left\|\mathcal{L}_{\mu_n}-\varphi_0(\mu_n)\right\|_V^2\right)\left(\sum_{n=1}^N\frac{{\left(g_{\mu},\Phi_n\right)^2_{\mathcal{U}}}}{\tau_n}\right)\label{eq:proofprop1a}\\
&\leq 4N\underset{\mu\in\mathcal{P}}{\sup}\left\|\mathcal{L}_{\mu}-\varphi_0(\mu)\right\|_V^2\left(\sum_{n=1}^N\frac{{\int_{\mu\in\mathcal{P}}\left(g_{\mu},\Phi_n\right)^2_{\mathcal{U}}}}{\tau_n}\right)\label{eq:proofprop1b}\\
&\leq 4|\mathcal{P}|N^2Z_N\left(\theta^{sN}_{\mathcal{L}}\right)^2,\label{eq:proofprop1c}
\end{align}
\end{subequations}
\end{itemize}
where we recall that $\displaystyle Z_N=\frac{1}{N}\sum_{n=1}^N\frac{{\int_{\mu\in\mathcal{P}}\left(g_{\mu},\Phi_n\right)^2_{\mathcal{U}}}}{\frac{|\mathcal{P}|}{N}\sum_{p=1}^N\left(g_{\mu_p},\Phi_n\right)^2_{\mathcal{U}}}$.
We recall that the considered POD has been carried-out on the discrete set $S=\left\{g_{\mu_l}\right\}_{1\leq l\leq N}$, see Table~\ref{tablePOD} for the definition and the properties of this POD.
Consider $\tilde{S}=\left\{g_{\mu}\right\}_{\mu\in\mathcal{P}}$: a second POD applied on this set leads to POD modes denoted 
$\tilde{\Phi}_n$ and eigenvalues given by
$\tilde{\tau}_n := \frac{1}{|\mathcal{P}|}\int_{\mu\in\mathcal{P}}\left(g_{\mu},\tilde{\Phi}_n\right)^2_{\mathcal{U}}$, $1\leq n\leq N$.
The POD decompositions on $S$ and $\tilde{S}$ are asymptotically equal when $N$ tends to infinity (this corresponds to the case where the dimension of subspace spanned by the elements of $\tilde{S}$ is infinite), leading to: $\forall n\geq 1$, $\Phi_n\underset{N\to\infty}{\rightarrow}\tilde{\Phi}_n$ and ${\tau}_n\underset{N\to\infty}{\rightarrow}\tilde{\tau}_n$.
Hence, we infer $\forall n\geq 1$, $\displaystyle\frac{{\int_{\mu\in\mathcal{P}}\left(g_{\mu},\Phi_n\right)^2_{\mathcal{U}}}}{\frac{|\mathcal{P}|}{N}\sum_{p=1}^N\left(g_{\mu_p},\Phi_n\right)^2_{\mathcal{U}}}\underset{N\to\infty}{\rightarrow}1$, leading to $Z_N\underset{N\to\infty}{\rightarrow}1$, then $Z_N$ is bounded.

We now conclude the proof from the control of the terms in~\eqref{eq:propc} and by noting that $\delta_{\hat{Q}_{m,d}}=0$ owing to Property~\ref{interpEIM}, and $\tau_{\hat{Q}_{m,d}}>0$.
\end{proof}

\begin{remark}[Boundedness of $(Z_N)_N$]
Another way to control the bound of $(Z_N)_N$ is to recognize a Riemann sum in the denominator of the quotient in the definition of $Z_N$, for which the sampling $\mu_p$ is selected by an EIM on $g$. The problem of defining the best sample of points to construct an interpolation is complex and is in general not solved, but the sample provided by EIM is
competitive compared to situations where the best behavior is known, see the numerical illustrations in~\cite{Maday}. In pratice, in~\cite{Maday} and in Figure~\ref{fig:DOE}, we observe the points selected by the EIM to be distributed quite regularly in the parameter space (in particular, the EIM cannot select twice the same point). Construct a Voronoi tesselation of $\mathcal{P}$ from this set of points, and denote $v_N^p$ the volume of the cells. Denote $M_N:=\frac{N}{|\mathcal{P}|}\underset{1\leq p\leq N}{\sup}v_N^p$, since $\frac{|\mathcal{P}|}{N}$ is the mean volume of the cells, the assumption of regular distribution for the EIM points leads to $M_N$ is close to $1$ for all $N$.
Then, $\displaystyle Z_N\leq M_N\frac{1}{N}\sum_{n=1}^N \frac{{\int_{\mu\in\mathcal{P}}\left(g_{\mu},\Phi_n\right)^2_{\mathcal{U}}}}{\sum_{p=1}^Nv_p^N\left(g_{\mu_p},\Phi_n\right)^2_{\mathcal{U}}}$, where $\displaystyle\sum_{p=1}^Nv_p^N\left(g_{\mu_p},\Phi_n\right)^2_{\mathcal{U}}$ is a Riemann sum converging to the integral $\displaystyle\int_{\mu\in\mathcal{P}}\left(g_{\mu},\Phi_n\right)^2_{\mathcal{U}}$ as $N$ tends to infinity.
\end{remark}

\renewcommand*{\proofname}{\textbf{Proof of Proposition~\ref{convergence3}}}
\begin{proof}
Using the triangular inequality, there holds 
\begin{equation}\label{eq:prop7demo}
\int_{\mu\in\mathcal{P}}\|\mathcal{X}_{\mu}^{N}-\mathcal{L}_\mu\|_V^{2}\leq 3\int_{\mu\in\mathcal{P}}\left\|\mathcal{L}_{\mu}-X_{m} g_{\mu}\right\|^{2}_{V}+3\int_{\mu\in\mathcal{P}}\left\|\sum_{n=1}^{N}\lambda_n(\mu)\left(\mathcal{L}_{\mu_n}-X_{m} g_{\mu_n}\right)\right\|^{2}_{V}+3\int_{\mu\in\mathcal{P}}\left\|X_m\left(I^{N}(g_\mu)-g_\mu\right)\right\|_V^2,
\end{equation}
where we recall that $X_{m} g_{\mu}$ is the $m$-th term of the considered power algorithm at parameter value $\mu$, converging to the quantity of interest $\mathcal{L}_\mu$, and the EIM approximation of $g$ is $I^N(g)(\vec{k},\mu)=\sum_{n=1}^N\lambda_n(\mu)g(\vec{k},\mu_n)$.
The second term in the right-hand side can be controlled in the same fashion as in the proof of Proposition~\ref{convergence1}: replacing $\mathcal{L}_{\mu}$ with $\mathcal{L}_{\mu}-X_{m} g_{\mu}$, there holds, denoting $\overrightarrow{\mathcal{L}-X_{m}g}=\left(\mathcal{L}_{\mu_n}-X_{m} g_{\mu_n}\right)_{1\leq n\leq N}$,
\begin{subequations}
\begin{align}
3\int_{\mu\in\mathcal{P}}\left\|\sum_{n=1}^{N}\lambda_n(\mu)\left(\mathcal{L}_{\mu_n}-X_{m} g_{\mu_n}\right)\right\|_V^2&\leq 6\int_{\mu\in\mathcal{P}}\left\|\left(G^t\vec{\lambda}(\mu)-\vec{h}(\mu)\right)^t\cdot\left(G^{-1}\left(\overrightarrow{\mathcal{L}-X_{m}g}\right)\right)\right\|_V^2\\
&\qquad+6\int_{\mu\in\mathcal{P}}\left\|\vec{h}^t(\mu)\cdot\left(G^{-1}\left(\overrightarrow{\mathcal{L}-X_{m}g}\right)\right)\right\|_V^2\label{eq:prop8a}\\
&\leq 24C_1^2(m)\frac{N\delta_N^2}{\tau_N} + 6|\mathcal{P}|C_1^2(m)N^2{Z_N},\label{eq:prop8b}
\end{align}
\end{subequations}
where we recall that $C_1(m)$ is a bound for the approximation of the considered power matrix algorithm, $\delta_N$ is the EIM approximation error, $(Z_N)_N$ is a bounded sequence, and $\tau_N$ denotes the $N$-th eigenvalue of the considered POD, see Table~\ref{tablePOD}.
The control of the first term in the right-hand side of~\eqref{eq:prop8a} was obtained in~\eqref{eq:control6d} replacing $\vec{\mathcal{L}}$ by $\overrightarrow{\mathcal{L}-X_{m}g}$, and the control of the second term was obtained in~\eqref{eq:proofprop1c} replacing $\varphi_0$ by $X_mg$.
The proof is ended by noting that $\delta_{\hat{Q}_{m,d}}=0$ and $\forall\mu\in\mathcal{P},~I^{\hat{Q}_{m,d}}(g_\mu)=g_\mu$ owing to Property~\ref{interpEIM}.
\end{proof}

\subsection{Comments}
\label{comments}

We recall the practical results of Section~\ref{mainresults} stated in Corollary~\ref{coro1}: the nonintrusive approximations for the inverse matrix~\eqref{eq:iter_solution_2} and for the log-det~\eqref{eq:iter_solution_3} are convergent with respect to the number of evaluations $N^{\rm EIM}$ of the quantity of interest.

Consider now $F$ a compact subset of $\mathcal{C}^0(\mathcal{P},V)$ containing $\mathcal{L}$. The second bound in Proposition~\ref{convergence1} can be weakened to the form: for all integer $m$,
\begin{equation}
\frac{1}{|\mathcal{P}|}\int_{\mu\in\mathcal{P}}\|\mathcal{X}_{\mu}^{\hat{Q}_{m,d}}-\mathcal{L}_\mu\|^{2}_V\leq 
4\left(1+{\hat{Q}_{m,d}}^2{Z_{\hat{Q}_{m,d}}}\right)\eta_{\mathcal{S}^{s\hat{Q}_{m,d}}}^2,
\end{equation}
where 
\begin{equation}
\eta_{\mathcal{S}^{s\hat{Q}_{m,d}}}:=\underset{v\in F}{\sup}~\underset{\varphi\in\mathcal{S}^{s{\hat{Q}_{m,d}}}}{\inf}~{\left\|v-\varphi\right\|_{\mathcal{C}^0(\mathcal{P},V)}},
\end{equation}
which is related to the following Kolmogorov $s{\hat{Q}_{m,d}}$-width:
\begin{equation}\label{nwidth}
d_{s{\hat{Q}_{m,d}}}(F,\mathcal{C}^0(\mathcal{P},V))=\underset{F^{s{\hat{Q}_{m,d}}}\subset\mathcal{C}^0(\mathcal{P},V)}{\rm inf}\quad\underset{v\in F}{\rm sup}\quad\underset{\varphi\in F^{s{\hat{Q}_{m,d}}}}{\rm inf}\quad\|v-\varphi\|_{\mathcal{C}^0(\mathcal{P},V)}=\underset{F^{s{\hat{Q}_{m,d}}}\subset\mathcal{C}^0(\mathcal{P},V)}{\rm inf}\quad\eta_{F^{s{\hat{Q}_{m,d}}}}.
\end{equation}
For the need of the proof, we considered the snapshots POD on the set $\left(g_{\mu_l}\right)_{1 \leq l \leq {\hat{Q}_{m,d}}}$, with values $\mu_n$ selected by a first EIM on $g$, which lead to a fixed projector $(\mathcal{J}^{\hat{Q}_{m,d}}v)(\mu)=\sum_{n=1}^{\hat{Q}_{m,d}} {\left(g_{\mu},\Phi_n\right)_{\mathcal{U}}}\sum_{p=1}^{\hat{Q}_{m,d}}G^{-1}_{np}v(\mu_p)$ for $v\in\mathcal{C}^0(\mathcal{P},V)$, and therefore a fixed subspace $\mathcal{S}^{s{\hat{Q}_{m,d}}}\subset\mathcal{C}^0(\mathcal{P},V)$, instead of the optimal subspace $F^{s{\hat{Q}_{m,d}}}$ in~\eqref{nwidth}.

In~\cite{convergenceGEIM}, upper bounds for the EIM error have been derived, for polynomial and exponential decay rates of the Kolmogorov $n-$width $d_n(\{g_\mu,\mu\in\mathcal{P}\}, \mathcal{U})$. In Proposition~\ref{convergence1} are made explicit the dependences on the EIM upper bound $\delta_{N}$ and on $\theta_{\mathcal{L}}^{sN}$, which is related to the approximation of $\mu\mapsto\mathcal{L}_\mu$ in $\mathcal{C}^0(\mathcal{P},V)$, not to the EIM approximation of $g_\mu$ in $\mathcal{U}$.

The convergence of the upper bounds in Proposition~\ref{convergence1} are difficult to observe in practice, due to difficulty of the numerical estimation of $\delta_{N}$ and $\theta_{\mathcal{L}}^{sN}$.
However, the convergence of the upper bound in~\eqref{eq:prop6b} seams reasonnable since, in reducible cases, the convergence of $N\theta_{\mathcal{L}}^{sN}$ is a mild assumption when $\theta_{\mathcal{L}}^{sN}$ is replaced by the Kolmogorov $sN$-width $d_{sN}(F,\mathcal{C}^0(\mathcal{P},V))$. 
In our numerical experiments, we observed that the EIM provides reasonnable approximation errors only in the case where $N=\hat{Q}_{m,d}={Q}_{m,d}$, probably due to the particular form of $\Kappa_{m,d}$, which is a discrete set of multi-indices, not a discrete sampling of some continuous variable: the elements $\vec{k}$ of $\Kappa_{m,d}$ seem to generate linearly independant elements $\mu\mapsto g_\mu(\vec{k})$ of $\mathcal{C}^0(\mathcal{P},V)$. This could be compared to the discrete EIM (DEIM), where the EIM algorithm is applied on the indices list of some POD vectors, which are all kept for the approximation~\cite{chaturantabut2010nonlinear}.
However, the main advantage of the EIM in our case is in the selection of relevant parameter values in a large set $\mathcal{P}_{\rm sample}$, see Figure~\ref{fig:DOE} showing the location of parameter values selected by EIM. In the numerical experiments of Section~\ref{sec:numerics}, we took $N={Q}_{m,d}$ for this reason, and we can assess the convergence of the approximation with respect to the upper bounds of Corollary~\ref{coro1}.

Nevertheless, we derived a rather general result of convergence in Proposition~\ref{convergence1}, which could be very usefull for other classes of models and problems. 
Notice also that with $\hat{Q}_{m,d}$ evaluation of the reference quantity $\mathcal{L}_{\mu}$, $\theta_{\mathcal{L}}^{s\hat{Q}_{m,d}}$ involves an approximation on a $s\hat{Q}_{m,d}$-dimensional subspace of $\mathcal{C}^0(\mathcal{P},V)$, where $s$ is the dimension of $V$.

\section{Numerical experiments}
\label{sec:numerics}

In this section, numerical comparisons between the presented algorithms and others taken from the literature are presented.

\subsection{Inverse operators and solution to linear systems}
\label{sec:numsollinsys}

\begin{figure}[h!] 
   \begin{center}
   \includegraphics[width=0.6\textwidth]{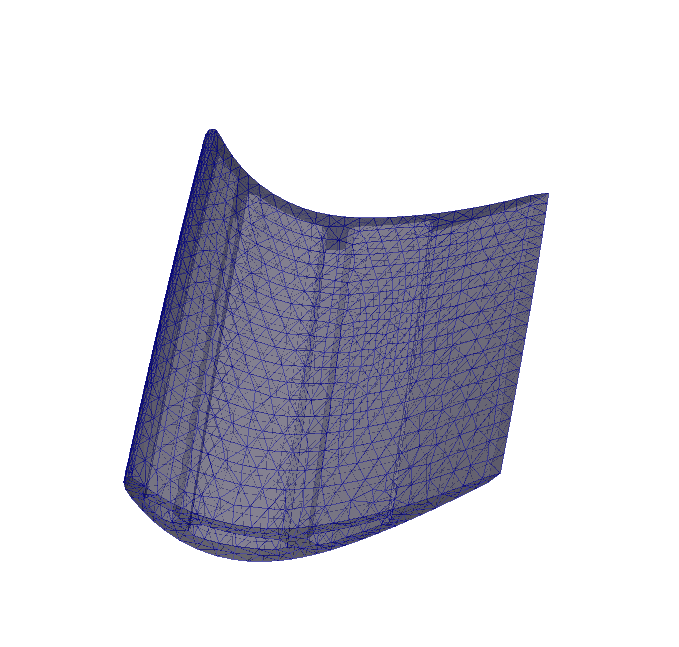}
   \end{center}
   \caption{\label{fig:mesh}Mesh of the high pressure turbine blade}
\end{figure}

Consider an open set $\Omega\subset\mathbb{R}^3$ meshed with tetrahedra, see Figure~\ref{fig:mesh}. 
This set represents a high pressure turbine blade featuring three cooling corridors; the intersection between these corridors and $\Omega$ is denoted $\partial\Omega_\mathcal{C}$. We consider the following archetypal heat problem:
\begin{equation}
\label{eq:formeforte}
\left\{
\begin{alignedat}{3}
-\vec{\nabla}\cdot\vec{q}&=\mu_2u_{\mu}\quad&&\textnormal{ in }\Omega,\\
\vec{q}\cdot\vec{n}&=-1\quad&&\textnormal{ on }\partial\Omega_{\mathcal{C}},&\\
\vec{q}\cdot\vec{n}&=0\quad&&\textnormal{ on }\partial\Omega\backslash\partial\Omega_{\mathcal{C}},&\\
\end{alignedat}
\right.
\end{equation}
where $u_{\mu}$ is the temperature, $\vec{q}=-\mu_1\vec{\nabla}u_{\mu}$ is the heat flux density, and $\mu=(\mu_1,\mu_2)\in\mathcal{P}:=(1,4)^2$ is the parameter. In this problem, $\mu_1$ is the heat conductivity, and $\mu_2 u_{\mu}$ is a volumic source term depending on the solution $u_{\mu}$.

\begin{figure}[H] 
   \begin{center}
   \includegraphics[width=\textwidth]{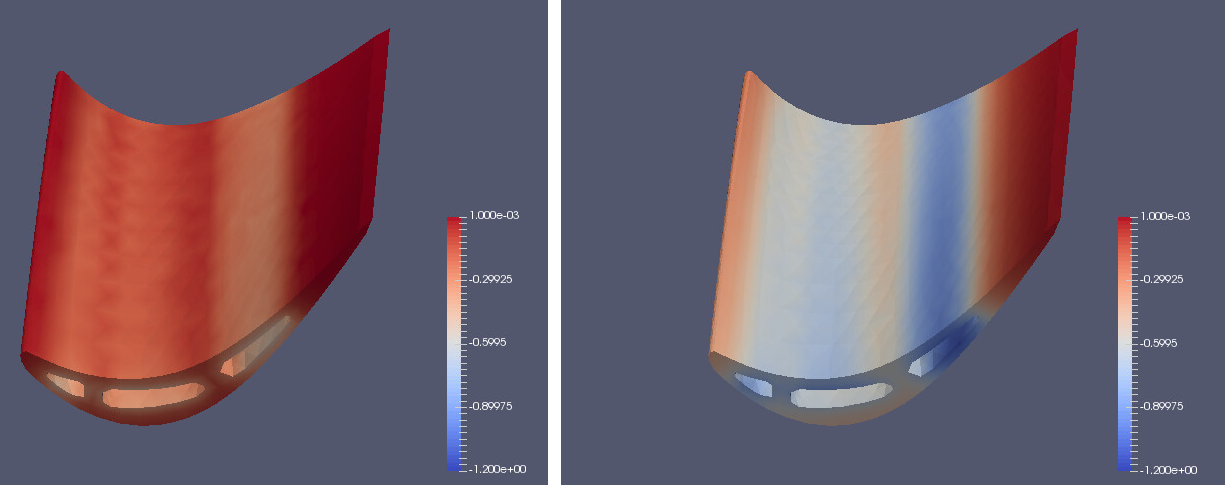}
   \end{center}
   \caption{\label{fig:solutions}Solutions $u_h$ to~\eqref{eq:formefaible}, for respective parameter values (1.82, 3.87) and (3.48,1.21)}
\end{figure}

Denote $\mathcal{V}_h(\Omega)$ the space of P1-finite elements associated with the considered mesh of $\Omega$, where $h$ denotes the characteristics length of the tetrahedra constituting the mesh. The weak form of~\eqref{eq:formeforte} can be approximated by
\begin{equation}
\label{eq:formefaible}
A_\mu U_{\mu} = b,
\end{equation}
where $A_\mu=\mu_1A_1+\mu_2A_2$, with $(A_1)_{i,j}=\int_\Omega \vec{\nabla}\phi_i\cdot\vec{\nabla}\phi_j$ and $(A_2)_{i,j}=\int_\Omega \phi_i\phi_j$, and $b_i = \int_{\partial\Omega_{\mathcal{C}}} \phi_i$; $\{\phi_i\}_{1\leq i\leq \mathcal{N}}$ denoting the P1-finite elements basis, where $\mathcal{N}=3,296$ in this example. The approximation ${u_{\mu}}_h\in\mathcal{V}_h$ of the solution $u_{\mu}$ of~\eqref{eq:formeforte} is obtained as ${u_{\mu}}_h:=\sum_{i=1}^\mathcal{N}{U_{\mu}}_i\phi_i$. Two solutions ${u_{\mu}}_h$ at two different parameter values are shown in Figure~\ref{fig:solutions}.

In this section, we compare the approximation~\eqref{eq:intersol} with other methods for approximating parametrized solutions:
\begin{enumerate}
\item(Minimisation in the Frobenius norm) Let $Y_{Q_{m,d}}={\rm Span}\{A_{\mu_1}^{-1},\cdots,A_{\mu_{Q_{m,d}}}^{-1}\}$, and define
\begin{equation}
\label{eq:Frobenius}
P_{Q_{m,d}}(\mu):=\underset{P\in Y_{Q_{m,d}}}{\rm argmin}\|I-PA(\mu)\|_F,
\end{equation}
where $\|\cdot\|_F$ denotes the Frobenius norm.
From~\cite{NouyZahm}, $P_{Q_{m,d}}(\mu)=\sum_{l=1}^{Q_{m,d}}\lambda_l(\mu)A_{\mu_l}^{-1}$, with $\lambda(\mu)$ the solution of 
\begin{equation}
\label{eq:linsysM}
M(\mu)\lambda(\mu)=S(\mu),
\end{equation}
with $M_{i,j}(\mu)={\rm trace}(A^T_\mu A_{\mu_i}^{-T}A_{\mu_j}^{-1}A_\mu)$ and $S_i(\mu)={\rm trace}(A_{\mu_i}^{-1}A_\mu)$. $P_{Q_{m,d}}(\mu)$ is then the best approximation of $A_{\mu}^{-1}$ expressed as a linear combination of inverses, in a chosen distance. However, the construction of $M(\mu)$ and $S(\mu)$ requires the online constructions of $A_{\mu_i}^{-1}A_\mu$, $1\leq i\leq Q_{m,d}$, which is too computationally expensive -- the goal of~\cite{NouyZahm} is to propose computationaly effective approximations of~\eqref{eq:Frobenius}. Here, we make use of~\eqref{eq:decompA0} to write $M_{i,j}(\mu)=\sum_{l,m=1}^{}\alpha_l(\mu)\alpha_m(\mu){\rm trace}(A_l^TA_{\mu_i}^{-T}A_{\mu_j}^{-1}A_m)$ and $S_i(\mu)=\sum_{l=1}^{Q_{m,d}}\alpha_l(\mu){\rm trace}(A_{\mu_i}^{-1}A_l)$, where the matrices
$\left({\rm trace}(A_l^TA_{\mu_i}^{-T}A_{\mu_j}^{-1}A_m)\right)_{i,j}$, $1\leq l,m\leq Q_{m,d}$ and right-hand sides $\left({\rm trace}(A_{\mu_i}^{-1}A_l)\right)_{i}$, $1\leq l\leq Q_{m,d}$ can be precomputed -- which is still much more computationally demanding than the proposed algorithm.

Finally, the approximation of~\eqref{eq:formefaible} can be written ${u_{\mu}}_h(\mu)\approx\sum_{l=1}^{Q_{m,d}}\lambda_l(\mu){u_{\mu}}_h(\mu_l)$, with $\lambda(\mu)$ the solution of~\eqref{eq:linsysM}. Notice that this method is intrusive since it requires to access the matrices $A_{\mu}$, instead of only the solutions ${u_{\mu}}_h(\mu)$ for nonintrusive strategies.

\item (Proper Orthogonal Decomposition (POD)~\cite{sirovich}) First, we construct $M\in\mathbb{R}^{Q_{n,d}\times N}$ such that $M_{i,j}={u_{\mu}}_h(\mu_i)_j$. Then, we compute the singular value decomposition (SVD) of $M$: $M=WSV$, with $S\in\mathbb{R}^{Q_{m,d}\times \mathcal{N}}$ containing the singular values of $M$ on its diagonal and zero elsewhere, and $W\in\mathbb{R}^{Q_{m,d}\times Q_{m,d}}$ and $V\in\mathbb{R}^{\mathcal{N}\times \mathcal{N}}$ unitary matrices. The $\hat{\mathcal{N}}$ largest singular values are kept, following an accuracy criterion, and we denote $\hat{V}\in\mathbb{R}^{\hat{\mathcal{N}}\times \mathcal{N}}$ such that $\hat{V}_{i,j}={V}_{i,j}$, $1\leq i\leq \hat{\mathcal{N}}$, $1\leq j\leq \mathcal{N}$, and $\{\hat{v}_i\}_{1\leq i\leq \hat{\mathcal{N}}}$, such that $\left(\hat{v}_i\right)_j=\hat{V}_{i,j}$, $1\leq i\leq \hat{\mathcal{N}}$, $1\leq j\leq \mathcal{N}$. Finally, the approximation of ${u_{\mu}}_h(\mu)$ is computed as ${u_{\mu}}_h(\mu)\approx\sum_{i=1}^{\hat{\mathcal{N}}}\theta_i(\mu)\hat{v}_i$, where $\theta(\mu)$ solves $\hat{A}_\mu \theta(\mu) = \hat{b}$, with $\hat{A}_\mu=\alpha_1(\mu)\hat{V}A_1\hat{V}^T+\alpha_2(\mu)\hat{V}A_2\hat{V}^T$, and $\hat{b}=\hat{V}b$.

Notice that this method is intrusive since it requires to access the matrices $A_{\mu}$, and that the result is not a linear combinations of solutions, but a linear combinations of singular vectors of $M$.

\item (Meta-modelisation) This first step is the same as the previous item: the construction of a basis $\hat{v}_i$, $1\leq i\leq \hat{\mathcal{N}}$ using the SVD. To provide a fair comparison, we did not use the snapshots ${u_{\mu}}_h(\mu_i)$ for the $\mu_i$ selected by the EIM on $g(\vec{k},\mu)$, but we used a Latin Hypercube Sampling of same size, computed with the MaxProj algorithm (see~\cite{maxproj}) to select the $\mu_i$, which is recommended when using statistical meta-models. The obtained set of parameter values is called the Design Of Experiment (DOE). Then, we compute the coefficients of the snapshots on the constructed basis: $\alpha_{i,j}=({u_{\mu_i}}_h,\hat{v}_j)$, $1\leq i\leq Q_{m,d}$, $1\leq j\leq \hat{\mathcal{N}}$.
Finally, nonintrusive meta-models are constructed on the coefficients $\alpha_{i,j}$, for which the predictions $\theta_i(\mu)$ at new parameter value $\mu$ are used in the obtained approximation as ${u_{\mu}}_h(\mu)\approx\sum_{i=1}^{\hat{\mathcal{N}}}\theta_i(\mu)\hat{v}_i$.
The considered statistical methods are taken from the machine learning community: (i) Gaussian processes, (ii) gradient boosting regression, (iii) random forests and (iv) Bayesian Ridge regression; and computed using the python package scikit-learn, see~\cite{scikit-learn}.
\end{enumerate}

\begin{figure}[H]
   \begin{center}
   \includegraphics[width=0.8\textwidth]{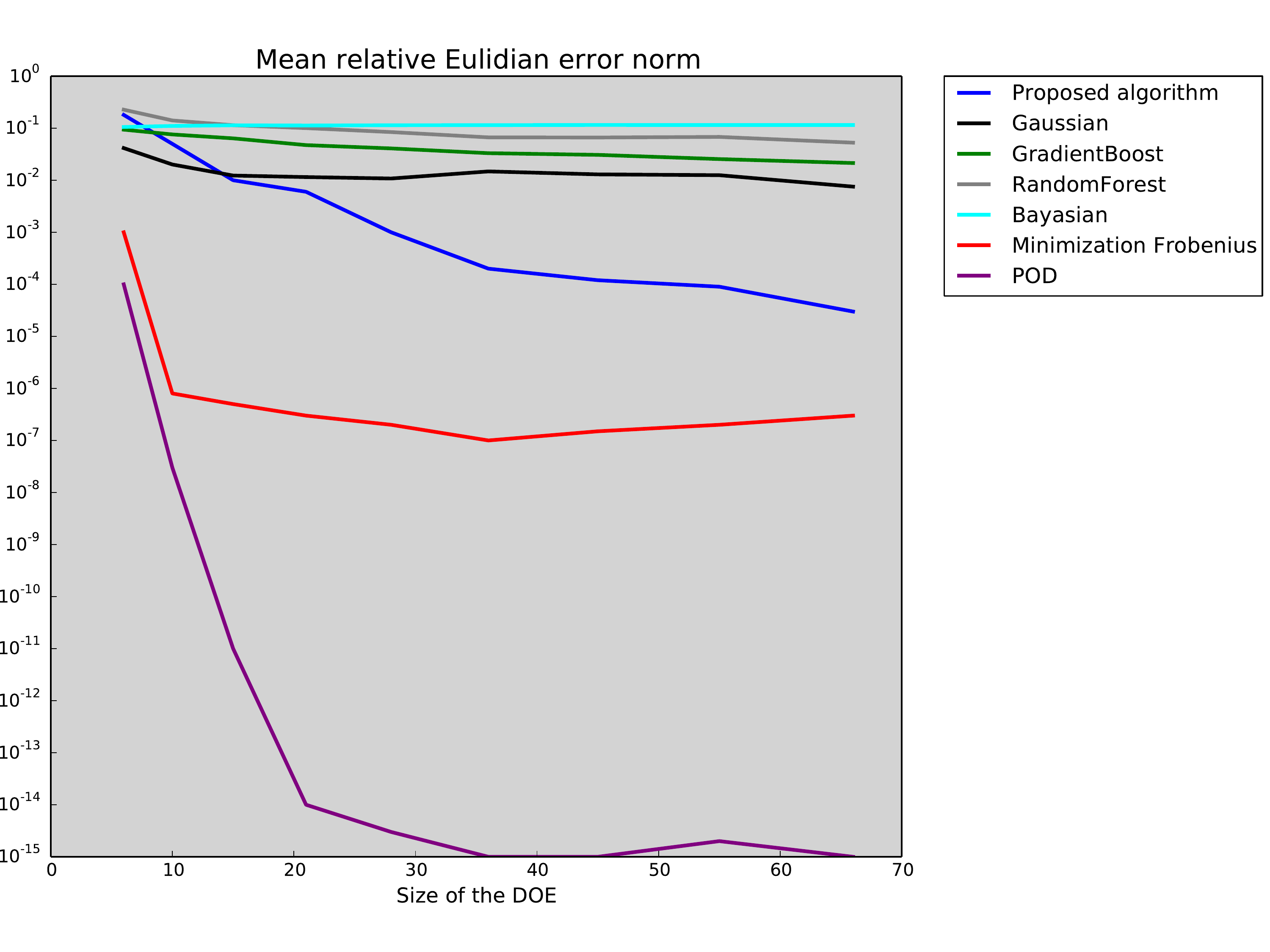}
   \end{center}
   \caption{\label{fig:fig1}Mean relative Euclidian-norms errors on the solution to~\eqref{eq:formefaible} using various interpolation methods, over a set of 100 randomly picked parameter values}
\end{figure}

Figure~\ref{fig:fig1} shows a comparison for the different aforementioned approximations. In abscissa is the size of the DOE for the statistical methods, and corresponds to $Q_{m,d}$, the number of parameter values selected by EIM for the proposed algorithm, the Minimization of Frobenius norm, and the POD.
We pick 100 random values for the parameter and compute the solutions $u_{\mu_t}$, $1\leq t\leq 100$. For each size of the DOE, we compare the predicted solution using the described approximations, at each value $\mu_t$, and compute the relative Euclidian-norms errors as $\displaystyle{\epsilon_{\mu_t}=\frac{\|u_{\mu_t}-\hat{u}_{\mu_t}\|_2}{\|u_{\mu_t}\|_2}}$, where $\hat{u}$ denotes here the considered approximation. Finally, the mean of relative errors, $\frac{1}{100}\sum_{t=1}^{100}\epsilon_{\mu_t}$, is represented on Figure~\ref{fig:fig1}.

We see that the intrusive methods that require to access the matrix $A_\mu$ perform much better than the nonintrusive methods in the example. Among the nonintrusive considered ones, our algorithm exhibits the best convergence rate.

\begin{figure}[H]
   \begin{center}
   \includegraphics[width=0.8\textwidth]{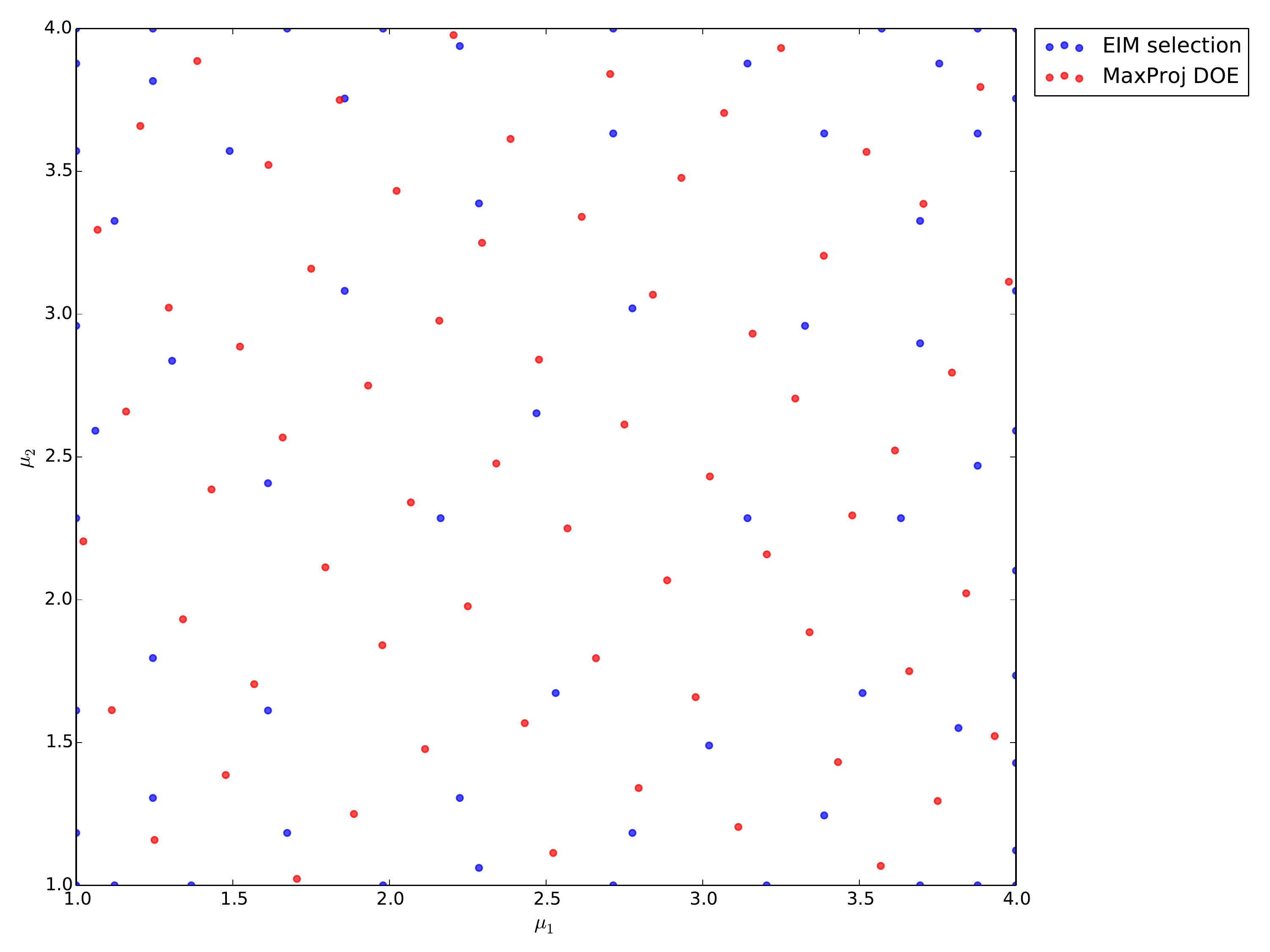}
   \end{center}
   \caption{\label{fig:DOE}Locations of the parameters selected by EIM and by the MaxProj DOE in the parameter domain $(1,4)^2$, for the largest DOE size considered in Figure~\ref{fig:fig1} (66 points)}
\end{figure}

Figure~\ref{fig:DOE} shows the locations of the parameters selected by EIM and by the MaxProj DOE in the parameter domain $(1,4)^2$. We notice that the EIM selects more points close to the boundary of the domain while the MaxProj DOE is more uniform.

\subsection{Logarithm of the determinant}

\begin{figure}[H] 
   \begin{center}
   \includegraphics[width=0.8\textwidth]{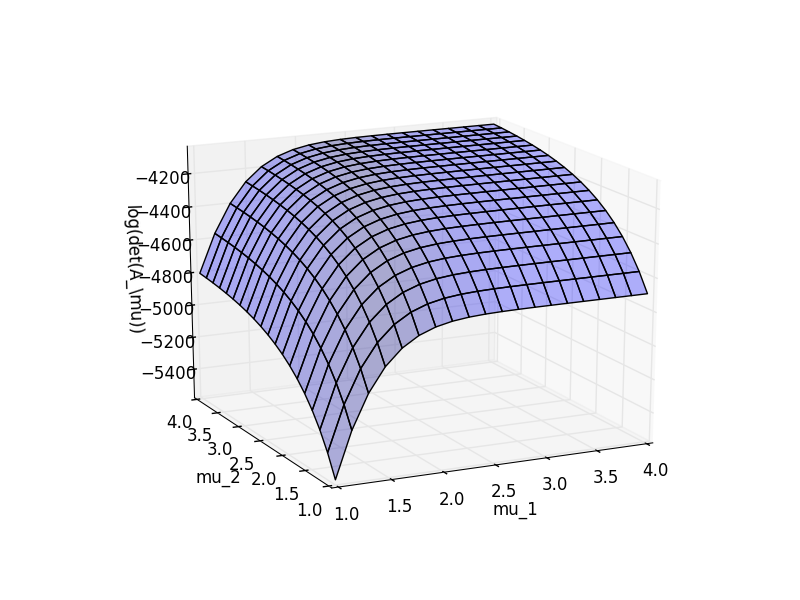}
   \end{center}
   \caption{\label{fig:logdetvalues}Representation of $\log(\det(A_\mu))$, with $\mu\in\mathcal{P}:=(1,4)^2$}
\end{figure}

In this section, we consider $A_\mu=0.045\left(1-e^{-\mu_1^2}\right)A_1+\left(1-e^{-\mu_2}\right)A_2$, where $A_1$ and $A_2$ are the same as in Section~\ref{sec:numsollinsys}, and $(\mu_1,\mu_2)\in\mathcal{P}:=(1,4)^2$. Figure~\ref{fig:logdetvalues} shows $\log(\det(A_\mu))$, $\mu\in\mathcal{P}$: this illustrates the quantity which we look to approximate using nonintrusive interpolation formulae (even if in large-dimensional test cases, we cannot afford to plot it).

\begin{figure}[H]
   \begin{center}
   \includegraphics[width=0.8\textwidth]{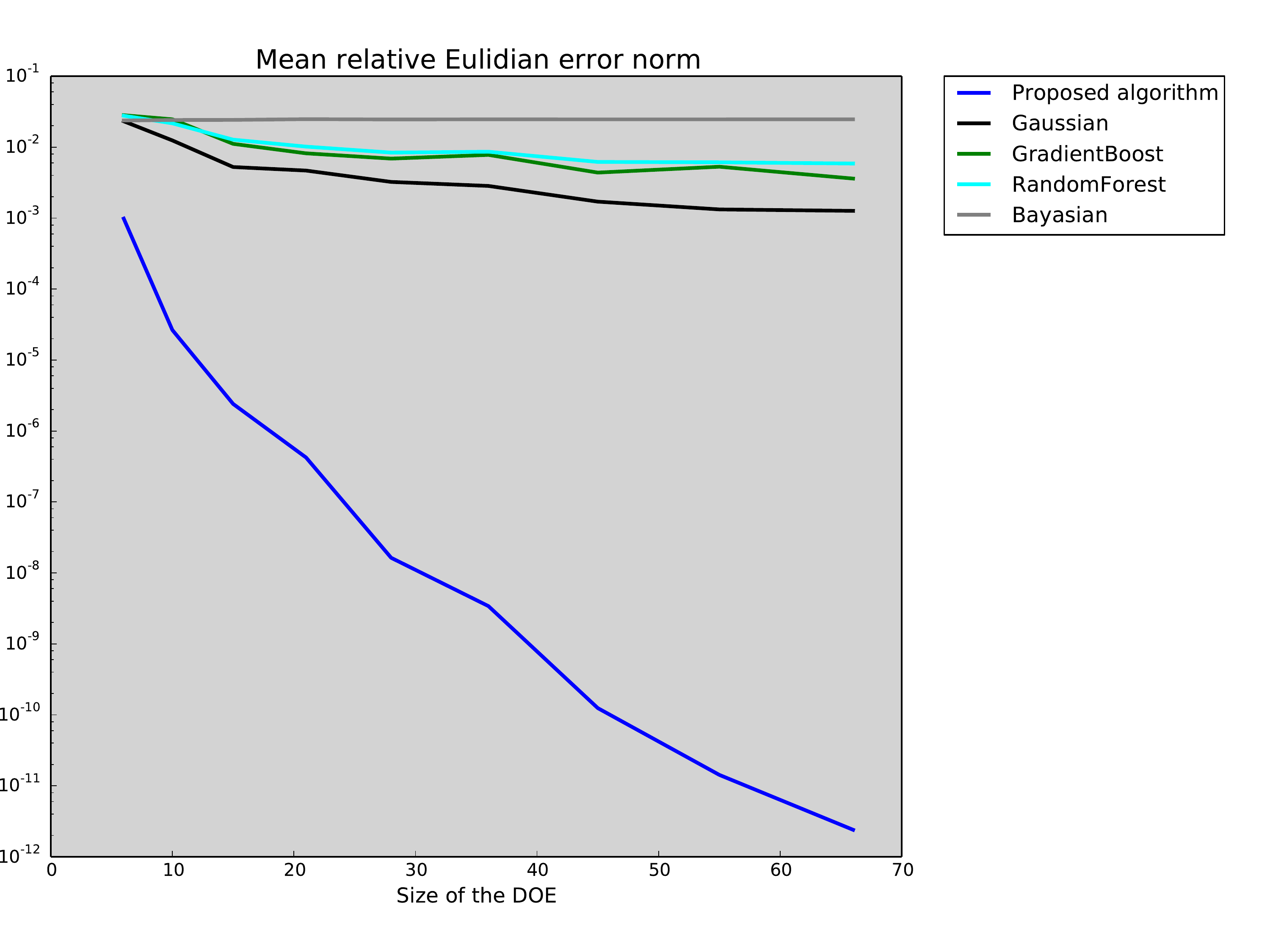}
   \end{center}
   \caption{\label{fig:logdetcomparision}Mean relative errors on the computation of $\log(\det(A_\mu))$ using the herein proposed algorithm and various nonintrusive machine learning regression methods}
\end{figure}

We compare the approximation~\eqref{eq:iter_solution_3} with the statistical nonintrusive approximation methods considered in Section~\ref{sec:numsollinsys}: (i) Gaussian processes, (ii) gradient boosting regression, (iii) random forests and (iv) Bayesian Ridge regression, see Figure~\ref{fig:logdetcomparision}. Here again, we construct a DOE using the MaxProj algorithm for selecting the values of the parameter used with the statistical nonintrusive approximation methods.
We pick 100 random values for the parameter and compute the solutions $u_{\mu_t}$, $1\leq t\leq 100$. For each size of the DOE, we compute $\displaystyle{\epsilon_{\mu_t}=\frac{\|\log(\det(A_{\mu_t}))-\hat{ld}_{\mu_t}\|_2}{\|\log(\det(A_{\mu_t}))\|_2}}$, where $\hat{ld}$ denotes here the considered approximation. Finally, the mean of relative errors, $\frac{1}{100}\sum_{t=1}^{100}\epsilon_{\mu_t}$, is represented on Figure~\ref{fig:logdetcomparision}.
We notice that the proposed algorithm exhibits the best convergence rate.

\subsection{Experiments in high parameter dimension cases}
\label{sec:highdimnum}

\subsubsection{A thermal problem in parameter dimension 10}
\label{sec:highdimtherm}

\begin{figure}[H] 
   \begin{center}
   \includegraphics[width=\textwidth]{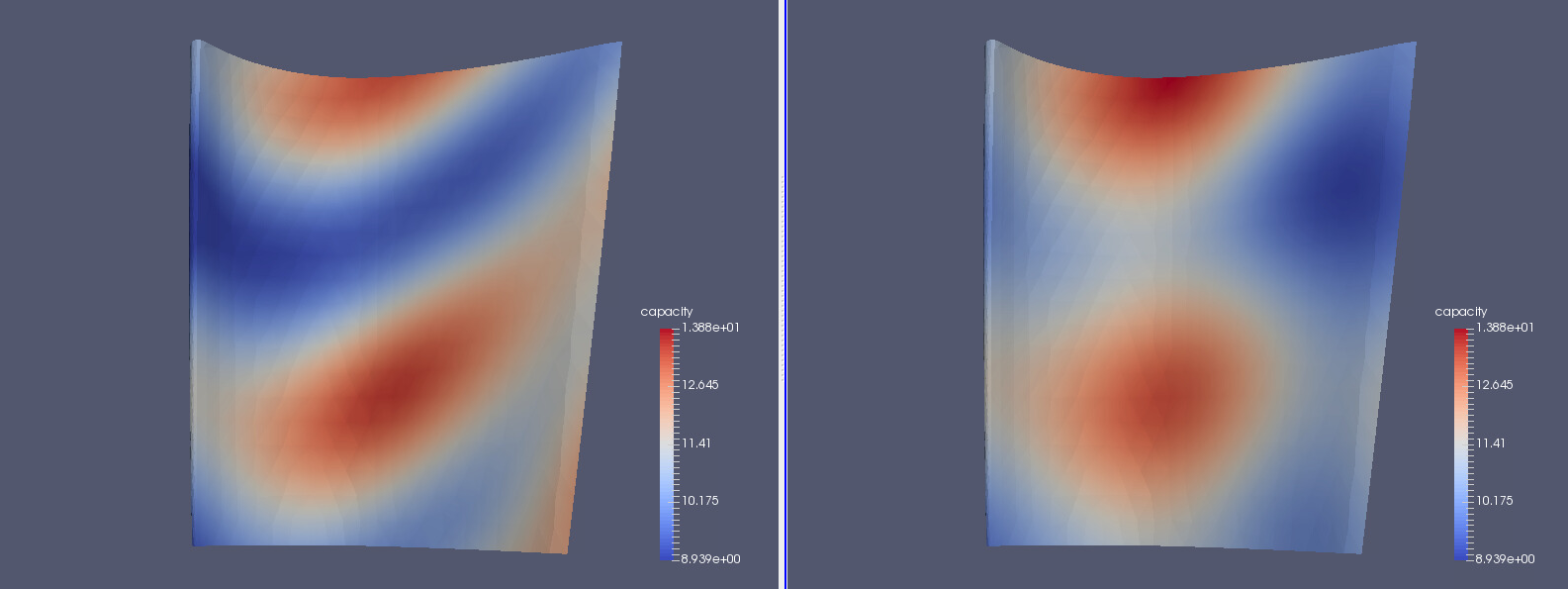}
   \end{center}
   \caption{\label{fig:2capacity}Example of heat capacities taken from experiment~1.}
\end{figure}

We consider the same geometry as in Section~\ref{sec:numsollinsys}, see Figure~\ref{fig:mesh}, and the following problem:
\begin{equation}
\label{eq:formeforte2}
\left\{
\begin{alignedat}{3}
c_p\frac{\partial u_{\mu}}{\partial t}-\vec{\nabla}\cdot\left(\eta\vec{\nabla}u_{\mu}\right)&=0\quad&&\textnormal{ in }\Omega,&\\
\eta\vec{\nabla}u_{\mu}&=1000~W.m^{-2}\quad&&\textnormal{ on }\partial\Omega_{\mathcal{C},}&\\
\vec{\nabla}u_{\mu}&=0\quad&&\textnormal{ on }\partial\Omega\backslash\partial\Omega_{\mathcal{C}},&\\
u_{\mu}&=0~{}^{\circ} C\quad&&\textnormal{ at }t=0,&
\end{alignedat}
\right.
\end{equation}
where $c_p$ denotes here the heat capacity multiplied by the density, $\eta$ is the thermal conductivity, and $u_{\mu}$ is the unknown temperature field. We choose $\eta=370~W.m^{-1}.K^{-1}$, and $c_p$ contains the parameter dependence as follows (in $J.m^{-3}.K^{-1}$):
\begin{itemize}
\item experiment~1:\quad $c_p(\mu,x)=10+\mu_1\cos(0.2x)+\mu_2\cos(0.25y)+\mu_3\cos(0.3z)+\mu_4\cos(0.2(x+y))+\mu_5\cos(0.25(x+z))+\mu_6\cos(0.3(y+z))+\mu_7\frac{x}{x_{\max}}+\mu_8\frac{y}{y_{\max}}+\mu_9\frac{z}{z_{\max}}+\mu_{10}\cos(0.1(x+y+z))$,\\
$\mu\in\mathcal{P}=(0.1,0.15)^{10}$, $(0.1,0.2)^{10}$, $(0.1,0.3)^{10}$, $(0.1,0.6)^{10}$ or $(0.1,1.1)^{10}$,
\item experiment~2:\quad$c_p(\mu,x)=10+(1-e^{-\mu_1})\cos(0.2x)+(1-e^{-\mu_2})\cos(0.25y)+(1-e^{-\mu_3})\cos(0.3z)+(1-e^{-\mu_4})\cos(0.2(x+y))+(1-e^{-\mu_5})\cos(0.25(x+z))+(1-e^{-\mu_6})\cos(0.3(y+z))+(1-e^{-\mu_7})\frac{x}{x_{\max}}+(1-e^{-\mu_8})\frac{y}{y_{\max}}+(1-e^{-\mu_9})\frac{z}{z_{\max}}+(1-e^{-\mu_{10}})\cos(0.1(x+y+z))$,\\
$\mu\in\mathcal{P}=(2,2.05)^{10}$, $(2,2.1)^{10}$, $(2,2.2)^{10}$, $(2,2.5)^{10}$ or $(2,3)^{10}$,
\end{itemize}
where $(x,y,z)\in\Omega$ is the space variable, and $x_{\max}$, $y_{\max}$ and $z_{\max}$ denote the maximum of the components of the space variable over $\Omega$. Examples of heat capacities taken from experiment~1 are represented in Figure~\ref{fig:2capacity}.

We are interested in the solution of~\eqref{eq:formeforte2} at $t=100 s$. Using a backward Euler time-discretization, the weak form reads: find $u_\mu\in H^1_0(\Omega)$ such that for all $v\in H^1_0(\Omega)$ ,
\begin{equation}
\int_{\Omega}c_p\left(\frac{u_\mu(\vec{x},t=100)-u_\mu(\vec{x},t=0)}{\Delta t}\right)v(\vec{x})+\int_{\Omega}\eta\vec{\nabla}u_\mu(\vec{x},t=100)\cdot\vec{\nabla}v(\vec{x})=
\int_{\partial\Omega_{\mathcal{C}}}\eta\vec{\nabla}u_\mu(\vec{x},t=100)\cdot\vec{n},
\end{equation}
where $\vec{n}$ is the exterior normal on $\partial\Omega_{\mathcal{C}}$, and $H^1_0(\Omega)=\{v\in L^2(\Omega)\textnormal{ such that }\vec{\nabla}v\in L^2(\Omega)\textnormal{ and }v|_{\partial\Omega\backslash\partial\Omega_{\mathcal{C}}}=0\}$. A finite element approximation is obtained as $A_\mu U_{\mu}=b$, where
\begin{equation}
\begin{alignedat}{3}
\left(A_\mu\right)_{i,j}&=\int_{\vec{x}\in\Omega}\frac{c_p(\vec{x},\mu)}{100}\phi_i(\vec{x})\phi_j(\vec{x})dx+370\int_{\vec{x}\in\Omega}\vec{\nabla}\phi_i(\vec{x})\cdot\vec{\nabla}\phi_j(\vec{x})d\vec{x},\qquad &1\leq i,j\leq\mathcal{N},\\
b_j&=1000\int_{\vec{x}\in\partial\Omega_{\mathcal{C}}}\phi_j(\vec{x})d\vec{x},\qquad &1\leq j\leq\mathcal{N},
\end{alignedat}
\end{equation}
where we recall that $\phi_i$ denoted the P1-finite element basis.
An approximation of~\eqref{eq:formeforte2} is obtained as $u_{\mu}(\vec{x},t=100)=\sum_{i=1}^{\mathcal{N}}{U_{\mu}}_i\phi_i(\vec{x})$. For experiment 1, the affine decomposition~\eqref{eq:decompA0} is obtained as $A_\mu=\sum_{l=1}^{11}\alpha_l(\mu)A_l$ with

\begin{small}
\hspace{-1.1cm}\begin{tabular}[b]{l l l}
   \multicolumn{2}{c}{$(A_1)_{i,j}=\int_{\vec{x}\in\Omega}\frac{1}{10}\phi_i(\vec{x})\phi_j(\vec{x})d\vec{x}+370\int_{\vec{x}\in\Omega}\vec{\nabla}\phi_i(\vec{x})\cdot\vec{\nabla}\phi_j(\vec{x})dx$} & $(A_2)_{i,j}=\int_{\vec{x}\in\Omega}\frac{\cos(0.2x)}{100}\phi_i(\vec{x})\phi_j(\vec{x})d\vec{x}$ \\
    $(A_3)_{i,j}=\int_{\vec{x}\in\Omega}\frac{\cos(0.25y)}{100}\phi_i(\vec{x})\phi_j(\vec{x})d\vec{x}$ & $(A_4)_{i,j}=\int_{\vec{x}\in\Omega}\frac{\cos(0.3z)}{100}\phi_i(\vec{x})\phi_j(\vec{x})d\vec{x}$ & $(A_5)_{i,j}=\int_{\vec{x}\in\Omega}\frac{\cos(0.2(x+y))}{100}\phi_i(\vec{x})\phi_j(\vec{x})d\vec{x}$ \\
    $(A_6)_{i,j}=\int_{\vec{x}\in\Omega}\frac{\cos(0.25(x+z))}{100}\phi_i(\vec{x})\phi_j(\vec{x})d\vec{x}$ & $(A_7)_{i,j}=\int_{\vec{x}\in\Omega}\frac{\cos(0.3(y+z))}{100}\phi_i(\vec{x})\phi_j(\vec{x})d\vec{x}$ & $(A_8)_{i,j}=\int_{\vec{x}\in\Omega}\frac{x}{100x_{\rm max}}\phi_i(\vec{x})\phi_j(\vec{x})d\vec{x}$ \\
    $(A_9)_{i,j}=\int_{\vec{x}\in\Omega}\frac{y}{100y_{\rm max}}\phi_i(\vec{x})\phi_j(\vec{x})d\vec{x}$ & $(A_{10})_{i,j}=\int_{\vec{x}\in\Omega}\frac{z}{100z_{\rm max}}\phi_i(\vec{x})\phi_j(\vec{x})d\vec{x}$ & $(A_{11})_{i,j}=\int_{\vec{x}\in\Omega}\frac{\cos(0.1(x+y+z))}{100}\phi_i(\vec{x})\phi_j(\vec{x})d\vec{x}$  
\end{tabular}
\end{small}
and $\alpha_1(\mu)=1$, $\alpha_{l+1}(\mu)=\mu_{(l)}$, $1\leq l\leq 10$ (notice that $\mu_{(l)}$ denotes the component $l$ for the continuous parameter $\mu\in\mathbb{R}^{10}$). For experiment 2, the matrices $A_l$, $1\leq l\leq 11$, are the same and $\alpha_1(\mu)=1$, $\alpha_{l+1}(\mu)=\left(1-e^{-\mu_{(l)}}\right)$, $1\leq l\leq 10$.

Figure~\ref{expTherm} shows the relative errors between the proposed algorithm and Gaussian processes, as detailed in Section~\ref{sec:numsollinsys}, in $L^2$- and $L^\infty$- norms. The DOE for the Gaussian processes in the parameter spaces is obtained using MaxProj as well, and we do not provide a comparison with the other statistical methods considered in the previous sections since they exhibited worse results. We notice that the proposed algorithm provides accurate results. Then, smaller parameter discrepancies lead to more accurate results: the reducibility of the problem is better, and $\rho$ in~\eqref{eq:rho} should be smaller leading to smaller bounds $C_1(m)$ in~\eqref{eq:coro1} (at each $m$). Moreover, at fixed parameter discrepancy (hence fixed~$\rho$), the errors decrease as $Q_{m,d}$ increases (hence as $m$ increases): the EIM computes exactly (i.e. with no approximation errors, due to the interpolation property~\ref{interpEIM}) more elements of the series defined by the iteration scheme~\eqref{eq:mat_it_scheme}, and the $A_{\mu_l}^{-1}$ in~\eqref{eq:iter_solution_2} are closer to the $X_m g_{\mu_l}$ in~\eqref{eq:iter_solution}. In~\eqref{eq:coro1}, this corresponds to the convergence of $C_1(m)$ to $0$ with respect to $m$.

\begin{figure}[h]
\begin{subfigure}{.44\linewidth}
    \centering
\resizebox{\linewidth}{\linewidth}{\input{dimension_10_testcase4_L2.tex}}
\caption{Relative $L^2$-norm error.}
    \end{subfigure}
    \qquad
\begin{subfigure}{.44\linewidth}
    \centering
\resizebox{\linewidth}{\linewidth}{\input{dimension_10_testcase4_Linf.tex}}
    \caption{Relative $L^\infty$-norm error.}
   \end{subfigure}

\vspace{0.5cm}
   
\begin{subfigure}{.44\linewidth}
    \centering
\resizebox{\linewidth}{\linewidth}{\input{dimension_10_testcase3_L2.tex}}
\caption{Relative $L^2$-norm error.}
    \end{subfigure}
    \qquad
\begin{subfigure}{.44\linewidth}
    \centering
\resizebox{\linewidth}{\linewidth}{\input{dimension_10_testcase3_Linf.tex}}
    \caption{Relative $L^\infty$-norm error.}
   \end{subfigure}
\caption{[(A-B): experiment 1 ; (C-D): experiment 2] Comparison between the proposed algorithm and Gaussian processes in a parameter dimension 10 test case: \protect\tikz{\protect\draw[red] (0,0)--(0.5cm,0);\protect\draw[white] (0,-0.1cm)--(0.5cm,-0.1cm);}~Gaussian 0.05, \protect\tikz{\protect\draw[blue] (0,0)--(0.5cm,0);\protect\draw[white] (0,-0.1cm)--(0.5cm,-0.1cm);}~Gaussian 0.1, \protect\tikz{\protect\draw[black] (0,0)--(0.5cm,0);\protect\draw[white] (0,-0.1cm)--(0.5cm,-0.1cm);}~Gaussian 0.2,\protect\tikz{\protect\draw[yellow] (0,0)--(0.5cm,0);\protect\draw[white] (0,-0.1cm)--(0.5cm,-0.1cm);}~Gaussian 0.5 ,\protect\tikz{\protect\draw[green] (0,0)--(0.5cm,0);\protect\draw[white] (0,-0.1cm)--(0.5cm,-0.1cm);}~Gaussian 1, \protect\tikz{\protect\draw[dashed,color=red] (0,0)--(0.5cm,0);\protect\draw[dashed,color=white] (0,-0.1cm)--(0.5cm,-0.1cm);}~Proposed algorithm 0.05, \protect\tikz{\protect\draw[dashed,color=blue] (0,0)--(0.5cm,0);\protect\draw[dashed,color=white] (0,-0.1cm)--(0.5cm,-0.1cm);}~Proposed algorithm 0.1, \protect\tikz{\protect\draw[dashed,color=black] (0,0)--(0.5cm,0);\protect\draw[dashed,color=white] (0,-0.1cm)--(0.5cm,-0.1cm);}~Proposed algorithm 0.2, \protect\tikz{\protect\draw[dashed,color=yellow] (0,0)--(0.5cm,0);\protect\draw[dashed,color=white] (0,-0.1cm)--(0.5cm,-0.1cm);}~Proposed algorithm 0.5, \protect\tikz{\protect\draw[dashed,color=green] (0,0)--(0.5cm,0);\protect\draw[dashed,color=white] (0,-0.1cm)--(0.5cm,-0.1cm);}~Proposed algorithm 1, the number indicates the parameter discrepancy.}\label{expTherm}
\end{figure}

\subsubsection{A mechanical problem in parameter dimension 14}

\begin{figure}[h!] 
   \begin{center}
   \includegraphics[width=0.8\textwidth]{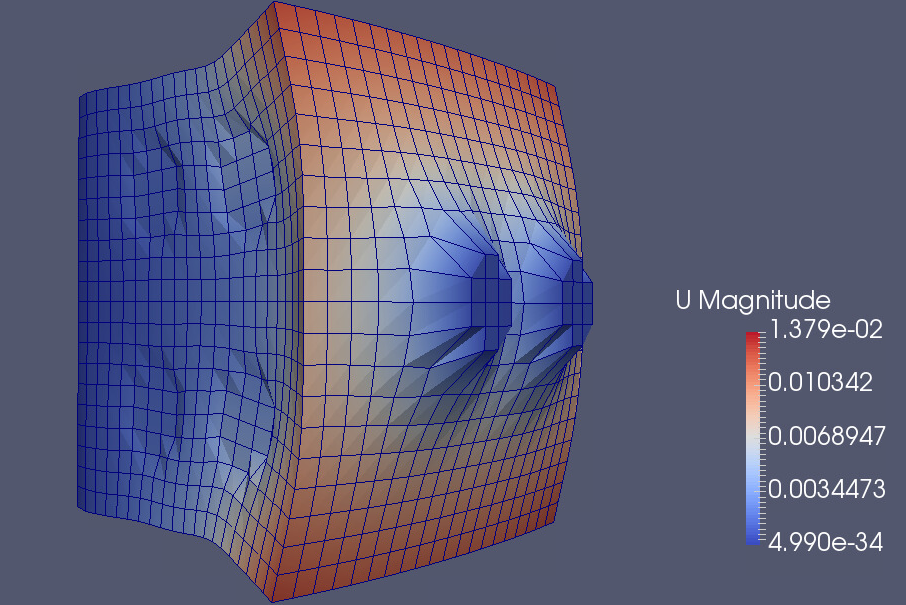}
   \end{center}
   \caption{\label{fig:mecamesh}Finite element approximation of the solution of Equation~\eqref{eq:variaformmeca}, where the mesh is deformed proportionally to the solution.}
\end{figure}

Consider a cube $\Omega$ meshed with linear hexahedra, with all displacement boundary conditions fixed on one face (denoted $\Gamma_D$) and a prescribed stress on the opposite face (denoted $\Gamma_N$), the other faces are free. The domain contains $6$ fibers $\Omega_1,\cdots,\Omega_6$, see Figure~\ref{fig:mecamesh}. We define $\Omega_0:=\Omega\backslash\left(\cup_{i=1}^6\Omega_i\right)$. We consider the following linear elasticity problem: find $u\in H^1_0(\Omega)^3$ such that $\forall v\in H^1_0(\Omega)^3$
\begin{equation}
\label{eq:variaformmeca}
\int_{\Omega}\frac{\eta_1}{2}\left(\nabla u_{\mu}+{}^t\nabla u_{\mu}\right)\cdot\left(\nabla v+{}^t\nabla v\right)+\int_{\Omega}\eta_2\left(\nabla\cdot u_{\mu}\right) \left(\nabla \cdot v\right)=\int_{\Gamma_N}(t\cdot n) v,
\end{equation}
where, $H^1_0(\Omega)^3=\{w\in L^2(\Omega)^3 \textnormal{ such that } \nabla w\in L^2(\Omega)^{3\times 3}\textnormal{ and } w|_{\Gamma_D}=0\}$, $\eta_1$ and $\eta_2$ are respectively Lam\'e's first and second parameters, $t=t_0n$ (with $n$ the outward unit normal and $t_0=-100~ N.m^{-2}$) is the prescribed traction vector on $\Gamma_N$, and ${u_{\mu}}$ is the unknown displacement. See Figure~\ref{fig:mecamesh} for a representation of a finite element approximation of the solution of Equation~\eqref{eq:variaformmeca}. We denote $\eta_{1,k}$ and $\eta_{2,k}$ respectively Lam\'e's first and second parameters of the subdomains $\Omega_k$, $0\leq k\leq 6$. We choose as parameter $\mu=\left(\eta_{1,0},\eta_{2,0},\eta_{1,1},\eta_{2,1},\cdots,\eta_{1,6},\eta_{2,6}\right)$, and the affine decomposition~\eqref{eq:decompA0} 
is obtained as $A_\mu=\sum_{l=1}^{14}\alpha_l(\mu)A_l$ with $(A_{2k})_{i,j}=\int_{\Omega_k}\frac{1}{2}\left(\nabla \phi_i+{}^t\nabla \phi_i\right)\cdot\left(\nabla \phi_j+{}^t\nabla \phi_j\right)$ and $(A_{2k+1})_{i,j}=\int_{\Omega_k}\left(\nabla\cdot \phi_i\right)\left(\nabla \cdot \phi_j\right)$, $0\leq k\leq 6$, $1\leq i,j\leq \mathcal{N}$ (where $(\phi_i)_{1\leq i\leq \mathcal{N}}$ is the basis of a finite element space approximating $H^1_0(\Omega)^3$, with $\mathcal{N}=27,783$), and $\alpha_{2k}=\eta_{1,k}$, $\alpha_{2k+1}=\eta_{2,k}$, $0\leq k\leq 6$. The parameter set is defined as follows: the reference Poisson coefficient is 0.3 in the whole cube, and the Young modulus for the fibers is $2\times10^{9}$, and $2\times10^{6}$ in the rest of the domain. From these values, we compute the reference Lam\'e's coefficients $(\eta_{1,k},\eta_{2,k})=(1.15\times10^9, 7.7\times10^8)$ for the fibers (namely $1\leq k\leq 6$), and $(\eta_{1,0},\eta_{2,0})=(1.15\times10^6, 7.7\times10^5)$ for the rest of the domain. Three parameter sets are considered, constituted of the intervals centered at the reference parameter values previously defined, with length respectively 1\%, 5\% and 10\% of the corresponding reference value.

In Figure~\ref{expMeca} are represented the relative errors between the proposed algorithm and Gaussian processes, as done in Section~\ref{sec:highdimtherm}, in $L^2$- and $L^\infty$- norms. The same conclusions as Section~\ref{sec:highdimtherm} can be drawn.

\begin{figure}[h!]
\begin{subfigure}{.44\linewidth}
    \centering
\resizebox{\linewidth}{\linewidth}{\input{dimension_14_testcase3_L2.tex}}
\caption{Relative $L^2$-norm error.}
    \end{subfigure}
    \qquad
    \begin{subfigure}{.44\linewidth}
    \centering
\resizebox{\linewidth}{\linewidth}{\input{dimension_14_testcase3_Linf.tex}}
    \caption{Relative $L^\infty$-norm error.}
   \end{subfigure}
\caption{Comparison between the proposed algorithm and Gaussian processes in a parameter dimension 14 test case: \protect\tikz{\protect\draw[red] (0,0)--(0.5cm,0);\protect\draw[white] (0,-0.1cm)--(0.5cm,-0.1cm);}~Gaussian 1\%, \protect\tikz{\protect\draw[blue] (0,0)--(0.5cm,0);\protect\draw[white] (0,-0.1cm)--(0.5cm,-0.1cm);}~Gaussian 5\%, \protect\tikz{\protect\draw[black] (0,0)--(0.5cm,0);\protect\draw[white] (0,-0.1cm)--(0.5cm,-0.1cm);}~Gaussian 10\%, \protect\tikz{\protect\draw[dashed,color=red] (0,0)--(0.5cm,0);\protect\draw[dashed,color=white] (0,-0.1cm)--(0.5cm,-0.1cm);}~Proposed algorithm 1\%, \protect\tikz{\protect\draw[dashed,color=blue] (0,0)--(0.5cm,0);\protect\draw[dashed,color=white] (0,-0.1cm)--(0.5cm,-0.1cm);}~Proposed algorithm 5\%, \protect\tikz{\protect\draw[dashed,color=black] (0,0)--(0.5cm,0);\protect\draw[dashed,color=white] (0,-0.1cm)--(0.5cm,-0.1cm);}~Proposed algorithm 10\%, the \% indicates the parameter discrepancy.}\label{expMeca}
\end{figure}

\section*{Conclusion}
In this work, we propose an algorithm to approximate, in a nonintrusive fashion, the limits of parametrized series of linear operators with respect to a functional $g$, based on the EIM approximation of $g$. We derive upper bounds of the error made by the obtained algorithm. With a strong enough convergence of the considered series, we prove the convergence of our algorithm. This assumption is verified by the two application considered in this work: the inverse and the logarithm of the determinant of a family of parametrized matrices. The numerical simulations illustrate that, in the considered test cases, our algorithm performs well compared to classical nonintruive approximations taken from the machine learning community.

\section*{Acknowledgement}
The authors would like to thank the anonymous reviewers for their relevant remarks and suggestions leading to significant improvements of the present work.
The author would also like to thank Tonya Rose from Safran for reviewing the manuscript.

\bibliographystyle{plain}
\bibliography{biblio}

\end{document}

%% file: dimension_10_testcase4_L2.tex
\begin{tikzpicture}

\definecolor{color0}{rgb}{0.75,0.75,0}

\begin{axis}[
xlabel={$Q_{m,d}$},
ylabel={Relative $L^2$-norm error},
xmin=9.34639255537567, xmax=336.600456417866,
ymin=6.0629935131116e-09, ymax=0.014265028539937,
xmode=log,
ymode=log,
xtick={0.1,1,10,100,1000,10000},
xticklabels={,,${10^{1}}$,${10^{2}}$,,},
ytick={1e-10,1e-09,1e-08,1e-07,1e-06,1e-05,0.0001,0.001,0.01,0.1,1},
yticklabels={,,${10^{-8}}$,${10^{-7}}$,${10^{-6}}$,${10^{-5}}$,${10^{-4}}$,${10^{-3}}$,${10^{-2}}$,,},
tick align=outside,
tick pos=left,
x grid style={white!69.019607843137251!black},
y grid style={white!69.019607843137251!black}
]
\addplot [red, dashed, forget plot]
table {%
11 1.98898608527487e-05
66 6.10509999055085e-08
286 1.18114884111328e-08
};
\addplot [red, forget plot]
table {%
11 0.000440275296046366
66 0.000295631343480467
286 0.000284963487266993
};
\addplot [blue, dashed, forget plot]
table {%
11 4.19950991006167e-05
66 6.00103691581969e-07
286 1.75228814095342e-08
};
\addplot [blue, forget plot]
table {%
11 0.00087537098142955
66 0.00058820102423471
286 0.000564400212591123
};
\addplot [black, dashed, forget plot]
table {%
11 0.000116332510549512
66 2.64187996461911e-06
286 8.84825556580779e-08
};
\addplot [black, forget plot]
table {%
11 0.00172456214177351
66 0.00116547730983866
286 0.00110804830769686
};
\addplot [color0, dashed, forget plot]
table {%
11 0.000802123428102758
66 3.66248139639e-05
286 1.3072262863407e-05
};
\addplot [color0, forget plot]
table {%
11 0.00402115766378482
66 0.00283498639040399
286 0.00262541180400994
};
\addplot [green!50.0!black, dashed, forget plot]
table {%
11 0.00115838050887134
66 0.000529752275648597
286 0.000214529637030907
};
\addplot [green!50.0!black, forget plot]
table {%
11 0.00732242817259765
66 0.00539391775577103
286 0.00486615926964923
};

\end{axis}

\end{tikzpicture}

%% file: dimension_10_testcase4_Linf.tex
\begin{tikzpicture}

\definecolor{color0}{rgb}{0.75,0.75,0}

\begin{axis}[
xlabel={$Q_{m,d}$},
ylabel={Relative $L^\infty$-norm error},
xmin=9.34639255537567, xmax=336.600456417866,
ymin=1.37120048811839e-07, ymax=0.058085646609091,
xmode=log,
ymode=log,
xtick={0.1,1,10,100,1000,10000},
xticklabels={,,${10^{1}}$,${10^{2}}$,,},
ytick={1e-08,1e-07,1e-06,1e-05,0.0001,0.001,0.01,0.1,1},
yticklabels={,,${10^{-6}}$,${10^{-5}}$,${10^{-4}}$,${10^{-3}}$,${10^{-2}}$,,},
tick align=outside,
tick pos=left,
x grid style={white!69.019607843137251!black},
y grid style={white!69.019607843137251!black}
]
\addplot [red, dashed, forget plot]
table {%
11 0.000225883721926886
66 9.47688381840895e-07
286 2.4709970445175e-07
};
\addplot [red, forget plot]
table {%
11 0.0022125624785752
66 0.00124727284426392
286 0.00127739216886634
};
\addplot [blue, dashed, forget plot]
table {%
11 0.000385562935756098
66 6.57984063148974e-06
286 3.80377730296124e-07
};
\addplot [blue, forget plot]
table {%
11 0.00438689314126953
66 0.00253533033023784
286 0.00254116275685253
};
\addplot [black, dashed, forget plot]
table {%
11 0.000886605073243631
66 2.13004982232331e-05
286 1.1757335621808e-06
};
\addplot [black, forget plot]
table {%
11 0.00857999166684946
66 0.00518910360172147
286 0.00518197691928786
};
\addplot [color0, dashed, forget plot]
table {%
11 0.00639198274645027
66 0.00048495683329656
286 0.000164050964988855
};
\addplot [color0, forget plot]
table {%
11 0.0189604607454319
66 0.0132229646186841
286 0.0135694779741877
};
\addplot [green!50.0!black, dashed, forget plot]
table {%
11 0.00695109634252115
66 0.00631959764093558
286 0.0029209406614421
};
\addplot [green!50.0!black, forget plot]
table {%
11 0.0322327649722503
66 0.0269583397330991
286 0.0282158285169142
};

\end{axis}

\end{tikzpicture}

%% file: dimension_10_testcase3_L2.tex
\begin{tikzpicture}

\definecolor{color0}{rgb}{0.75,0.75,0}

\begin{axis}[
xlabel={$Q_{m,d}$},
ylabel={Relative $L^2$-norm error},
xmin=9.34639255537567, xmax=336.600456417866,
ymin=4.00524874438311e-09, ymax=0.00105206421653751,
xmode=log,
ymode=log,
xtick={0.1,1,10,100,1000,10000},
xticklabels={,,${10^{1}}$,${10^{2}}$,,},
ytick={1e-10,1e-09,1e-08,1e-07,1e-06,1e-05,0.0001,0.001,0.01,0.1},
yticklabels={,,${10^{-8}}$,${10^{-7}}$,${10^{-6}}$,${10^{-5}}$,${10^{-4}}$,${10^{-3}}$,,},
tick align=outside,
tick pos=left,
x grid style={white!69.019607843137251!black},
y grid style={white!69.019607843137251!black}
]
\addplot [red, dashed, forget plot]
table {%
11 2.2751163553454e-07
66 7.82301176431729e-09
286 2.7173410700975e-08
};
\addplot [red, forget plot]
table {%
11 4.55741346126495e-05
66 3.34109228586944e-05
286 3.02692640334513e-05
};
\addplot [blue, dashed, forget plot]
table {%
11 6.01099152714677e-07
66 3.68393865634896e-08
286 6.79007998418039e-08
};
\addplot [blue, forget plot]
table {%
11 8.93655770155118e-05
66 6.51665726531114e-05
286 5.90692793192622e-05
};
\addplot [black, dashed, forget plot]
table {%
11 8.737675487114e-07
66 7.06262970029431e-09
286 1.44173293876284e-08
};
\addplot [black, forget plot]
table {%
11 0.000169226261566155
66 0.0001239371402937
286 0.000112524776646743
};
\addplot [color0, dashed, forget plot]
table {%
11 6.09826887741581e-06
66 3.27322765041008e-08
286 1.92354815592204e-08
};
\addplot [color0, forget plot]
table {%
11 0.000367951216001953
66 0.000267654415317028
286 0.000243832044917355
};
\addplot [green!50.0!black, dashed, forget plot]
table {%
11 1.51080650977701e-05
66 3.87127659400942e-07
286 1.10136985822462e-08
};
\addplot [green!50.0!black, forget plot]
table {%
11 0.000596630300767668
66 0.000422697144546918
286 0.000390244403805179
};

\end{axis}

\end{tikzpicture}

%% file: dimension_10_testcase3_Linf.tex
\begin{tikzpicture}

\definecolor{color0}{rgb}{0.75,0.75,0}

\begin{axis}[
xlabel={$Q_{m,d}$},
ylabel={Relative $L^\infty$-norm error},
xmin=9.34639255537567, xmax=336.600456417866,
ymin=7.7385709064192e-08, ymax=0.00434868427370505,
xmode=log,
ymode=log,
xtick={0.1,1,10,100,1000,10000},
xticklabels={,,${10^{1}}$,${10^{2}}$,,},
ytick={1e-09,1e-08,1e-07,1e-06,1e-05,0.0001,0.001,0.01,0.1},
yticklabels={,,${10^{-7}}$,${10^{-6}}$,${10^{-5}}$,${10^{-4}}$,${10^{-3}}$,,},
tick align=outside,
tick pos=left,
x grid style={white!69.019607843137251!black},
y grid style={white!69.019607843137251!black}
]
\addplot [red, dashed, forget plot]
table {%
11 2.70292804790997e-06
66 1.27220204551472e-07
286 4.73478110028003e-07
};
\addplot [red, forget plot]
table {%
11 0.000220815314718456
66 0.000143891177184569
286 0.000157145708492195
};
\addplot [blue, dashed, forget plot]
table {%
11 5.83408862455479e-06
66 5.90532483901817e-07
286 1.20224152047394e-06
};
\addplot [blue, forget plot]
table {%
11 0.000427015548146997
66 0.000279664123321011
286 0.000305879589914927
};
\addplot [black, dashed, forget plot]
table {%
11 5.56824428833226e-06
66 1.4413367461421e-07
286 2.4934779574564e-07
};
\addplot [black, forget plot]
table {%
11 0.000791424142864886
66 0.000528907857411894
286 0.000580263347865485
};
\addplot [color0, dashed, forget plot]
table {%
11 5.39877047645195e-05
66 3.97250668792434e-07
286 5.26722899377067e-07
};
\addplot [color0, forget plot]
table {%
11 0.00165391639002776
66 0.00113689376296638
286 0.00124187561678088
};
\addplot [green!50.0!black, dashed, forget plot]
table {%
11 0.000112423392203977
66 5.20211361919861e-06
286 1.31816644380414e-07
};
\addplot [green!50.0!black, forget plot]
table {%
11 0.00264522461037869
66 0.00190138656171047
286 0.00200721430866267
};

\end{axis}

\end{tikzpicture}

%% file: dimension_14_testcase3_L2.tex
\begin{tikzpicture}

\begin{axis}[
xlabel={$Q_{m,d}$},
ylabel={Relative $L^2$-norm error},
xmin=12.3956805755994, xmax=822.867283308226,
ymin=8.33607367630474e-10, ymax=0.0782313051019273,
xmode=log,
ymode=log,
xtick={1,10,100,1000,10000},
xticklabels={,,${10^{2}}$,,},
ytick={1e-11,1e-10,1e-09,1e-08,1e-07,1e-06,1e-05,0.0001,0.001,0.01,0.1,1},
yticklabels={,,${10^{-9}}$,${10^{-8}}$,${10^{-7}}$,${10^{-6}}$,${10^{-5}}$,${10^{-4}}$,${10^{-3}}$,${10^{-2}}$,,},
tick align=outside,
tick pos=left,
x grid style={white!69.019607843137251!black},
y grid style={white!69.019607843137251!black},
]
\addplot [red, dashed]
table {%
15 0.00140920781427981
120 1.82896659716133e-06
680 1.92019010598944e-09
};
\addplot [red]
table {%
15 0.000626396679822821
120 0.000567616243989389
680 0.000586097685087872
};
\addplot [blue, dashed]
table {%
15 0.033962362429063
120 1.18914156385539e-05
680 9.81598095481165e-07
};
\addplot [blue]
table {%
15 0.00169826595512778
120 0.0011166657820412
680 0.00110624056044591
};
\addplot [black, dashed]
table {%
15 0.0159586601512966
120 8.38251717898582e-05
680 9.1580958685816e-06
};
\addplot [black]
table {%
15 0.00299790720692941
120 0.0018485493315393
680 0.0018419986712284
};
\end{axis}
\end{tikzpicture}

%% file: dimension_14_testcase3_Linf.tex
\begin{tikzpicture}

\begin{axis}[
xlabel={$Q_{m,d}$},
ylabel={Relative $L^\infty$-norm error},
xmin=12.3956805755994, xmax=822.867283308226,
ymin=3.38918705552854e-09, ymax=0.24739913654271,
xmode=log,
ymode=log,
xtick={1,10,100,1000,10000},
xticklabels={,,${10^{2}}$,,},
ytick={1e-10,1e-09,1e-08,1e-07,1e-06,1e-05,0.0001,0.001,0.01,0.1,1,10},
yticklabels={,,${10^{-8}}$,${10^{-7}}$,${10^{-6}}$,${10^{-5}}$,${10^{-4}}$,${10^{-3}}$,${10^{-2}}$,${10^{-1}}$,,},
tick align=outside,
tick pos=left,
x grid style={white!69.019607843137251!black},
y grid style={white!69.019607843137251!black},
]
\addplot [red, dashed]
table {%
15 0.0050073182612182
120 9.17306402322103e-06
680 7.71824190884197e-09
};
\addplot [red]
table {%
15 0.00187099821008751
120 0.00192920553032225
680 0.00199986840408929
};
\addplot [blue, dashed]
table {%
15 0.108636391683828
120 4.06796258777552e-05
680 4.23661558241123e-06
};
\addplot [blue]
table {%
15 0.00452977404337265
120 0.00410330086835397
680 0.00445173005425543
};
\addplot [black, dashed]
table {%
15 0.0455886318638798
120 0.000267731550757115
680 4.78551633722374e-05
};
\addplot [black]
table {%
15 0.00890568825818388
120 0.00813193051901475
680 0.00880940824637026
};

\end{axis}

\end{tikzpicture}